\documentclass[11pt,seceqn,secthm]{elsart}
\usepackage{cite}
\usepackage{amsmath}
\usepackage{amsfonts}
\usepackage{amssymb}
\usepackage{amscd}
\journal{Topology}
\input{epsf.sty}

 \theoremstyle{plain}
\newtheorem{defprop}[thm]{Definition-Proposition}

\newenvironment{df}
{\vskip 5pt \addtocounter{thm}{1} \noindent{\bf Definition \thethm}
}{ \vskip 5pt}

\newenvironment{rmk}
{\vskip 5pt \addtocounter{thm}{1} \noindent{\bf Remark \thethm} }{
\vskip 5pt}

\newenvironment{nota}
{\vskip 5pt \addtocounter{thm}{1} \noindent{\bf Notation \thethm} }{
\vskip 5pt}

\newenvironment{quest}
{\vskip 5pt \addtocounter{thm}{1} \noindent{\bf Question \thethm} }{
\vskip 5pt}

\newenvironment{proof}
{\vskip 5pt \noindent{\bf Proof.}}{\hfill $\square$ \vskip 5pt}

\def\nn{\nonumber}

\def\br{br}
\def\cut{cut}
\def\gr{gr}
\def\In{In}

\def\Sh{Sh}
\def\Sn{\mathbb{S}_n}

\def\Stau{\mathbb{S}(\tau)}

\def\Z{\mathbb{Z}}

\def\ra{\rightarrow}
\def\del{\partial}
\def\s{\sigma}
\def\D{\Delta}
\def\SE{S}

\def\iso{\simeq}
\def\t{\tau}
\def\G{\Gamma}

\def\lb{L}
\def\verttoedge{{out}}
\def\whitevert{v_{white}}

\def\num{{num}}
\def\pin{{pin}}

\def\plant{{plant}}

\def\bp{bp}
\def\sign{{sign}}
\def\forgetainfty{st_{\infty}}
\def\cppin{cppin}

\def\cell{C_{\infty}}

\def\color{{clr}}
\def\op{op}

\def\nat{{Nat}}
\def\bnat{\overline{{Nat}}}
\def\Op{{Op}}
\def\bop{\overline{Op}}
\def\lab{{Lab}}
\def\blab{\overline{Lab}}

\def\rtree{\mathcal{T}^r}

\def\bwpptree{\mathcal{T}_{b/w}^{pp}}
\def\bppptree{\mathcal{T}_{bp}^{pp}}
\def\bwbptree{\mathcal{T}_{bp}^{pp}}
\def\stpptree{\mathcal{T}_{st}^{pp}}
\def\flrtree{\mathcal{T}^{r,fl}}

\def\wptree{\mathcal{T}^{pp,nt,fl}_{b/w}}
\def\wpttree{\mathcal{T}^{pp,fl}_{b/w}}

\def\wlptreet{\mathcal{T}_{b/w}^{pp,st}}
\def\wlbptreet{\mathcal{T}_{bp}^{pp}}
\def\wlptree{\mathcal{T}_{b/w}^{pp,st,nt}}
\def\wlbptree{\mathcal{T}_{bp}^{pp,nt}}

\def\H{\mathcal{H}opf}

\def\Set{\mathcal{S}et}
\def\Top{\mathcal{T}op}
\def\Vect{\Vectk}
\def\Vectk{\mathcal{V}ect_k}
\def\Chain{\mathcal{C}hain}

\def\Hom{\mathcal{H}om}
\def\Arc{\mathcal{A}rc}
\def\O{\mathcal{O}}
\def\cact{\mathcal{C}act}
\def\brace{\mathcal{B}race}

\def\cacti{\mathcal{C}acti}
\def\Cact{\mathcal{C}act}
\def\darc{\mathcal{DA}rc}
\def\DArc{\mathcal{DA}rc}
\def\Loop{\mathcal{L}oop}

\def\Pcact{\mathbb{P}\mathcal{C}act}

\def\CWCact{K}
\def\CCcact{CC_*(\Cact^1)}
\def\CCcactn{CC_*(\Cact^1(n))}
\def\CCkcactn{CC_k(\Cact^1(n))}
\def\CCncactn{CC_n(\Cact^1(n))}

\begin{document}

\begin{frontmatter}
\title
{On Spineless Cacti, Deligne's Conjecture and Connes-Kreimer's Hopf
Algebra}

\author
{Ralph M.\ Kaufmann\thanksref{me}}
\thanks[me]{Partially supported by NSF grant DMS-\#0070681
}

\ead{kaufmann@math.uconn.edu}

\address{University of Connecticut, Department of Mathematics, 196 Auditorium Rd., Storrs, CT 06269 }

\begin{abstract}
Using a cell model for the little discs operad in terms of
spineless cacti we give a minimal common topological operadic
formalism for three { a priori} disparate algebraic structures:
(1) a solution to Deligne's conjecture on the Hochschild complex,
(2) the Hopf algebra of Connes and Kreimer, and (3) the string
topology of Chas and Sullivan.
\end{abstract}

\end{frontmatter}


\section*{Introduction}

When considering an algebraic structure there is often a topological
framework which is indicative of this structure. For instance, it is
well known that Gerstenhaber algebras are governed by the homology
operad of little discs operad \cite{cohen,cohen2} and that
Batalin-Vilkovisky algebras are exactly the algebras over the
homology of the framed little discs operad \cite{Ge}. The purpose of
this paper is to prove that there is a common, minimal, topological
operadic formalism \cite{KLP,cact} for three { a priori} disparate
algebraic structures: (1) a homotopy Gerstenhaber structure on  the
chains of the Hochschild complex of an associative algebra, a.k.a.\
Deligne's conjecture, (2) the Hopf algebra of Connes and Kreimer
\cite{CK}, and (3) string topology \cite{CS1}.

To accomplish this task, we use the spineless cacti operad of
\cite{cact} which is responsible for the Gerstenhaber structure of
string topology \cite{CS1,Vor,KLP,cact,cyclic}. Moreover our operad of
spineless cacti is actually equivalent to the little discs operad
\cite{cact}.  The analysis of this operad on the chain level
allows us to give a new topological proof of Deligne's conjecture.
Furthermore it provides chain models for the operads whose
algebras are precisely pre-Lie algebras and graded pre-Lie
algebras, respectively, as well as a chain realization for the Hopf
algebra of Connes and Kreimer.

As we are using cacti, this approach naturally lies within string
topology on one hand and on the other hand it  is embedded in the
framework of a combinatorial description of the moduli space of
surfaces with punctured boundaries via the $\Arc$ operad
\cite{cact,KLP,P2}. Therefore all the previous
structures obtain a representation in terms of moduli spaces.

We start by giving new CW decompositions for the spaces
$\Cact^1(n)$ of normalized spineless cacti with $n$ lobes which
are homotopy equivalent to the spaces $\Cact(n)$ of spineless
cacti with $n$ lobes, the homotopy being the contraction of $n$
factors of $\mathbb{R}_{>0}$.

\noindent {\bf Theorem \ref{cact1}} {\it  The space $\Cact^1(n)$ is
homeomorphic to
 the CW complex $\CWCact(n)$.}


As shown in \cite{cact} the spaces $\cact^1(n)$ form a quasi-operad
whose homology is an operad isomorphic to the homology operad of cacti
and hence to the homology of the little discs operad. The operad structure, however,
already appears on the chain level.

\noindent {\bf Theorem \ref{celldecomp}} {\it The glueings induced
from the glueings of spineless normalized cacti make the spaces
$\CCcactn$ into a chain operad. Thus $\CCcact$ is an operadic model
for the chains of the little discs operad.}


Moreover, the cells of $\CWCact(n)$ are indexed by planted planar
bipartite trees and the operad of cellular chains is isomorphic to
a combinatorial tree dg-operad. Reinterpreting the trees as ``flow
charts'' for multiplications and brace operations and specifying
appropriate signs,  we obtain an operation of the cell operad of
cacti and hence a cell model of the little discs operad on the
Hochschild cochains of an associative algebra. This proves
Deligne's conjecture in any characteristic, i.e.\ over
$\mathbb{Z}$.

\noindent {\bf Theorem \ref{delconj}} {\it Deligne's conjecture is
true for the chain model of the little discs operad provided by
$\CCcact$, that is $CH^*(A,A)$ is a dg-algebra over $\CCcact$
lifting the Gerstenhaber algebra structure. }


Moreover, all possible flow charts using multiplication and brace
operations are realized by the operations of the cells, and these
operations are exactly the set of operations which appear when
studying iterations of the bracket and the product on the
Hochschild cochains. In this sense our solution to Deligne's
conjecture is minimal.

Deligne's conjecture has by now been proven in various ways
\cite{Maxim,T,MS,Vor2,KS,MS2,BF} (for a full review of the history
see \cite{MSS}). The different approaches  are basically realized
by choosing adequate chain models and some more or less abstract
form of homological algebra. The virtue of our approach which is
in spirit close to those of \cite{MS} and \cite{KS} lies in its
naturality and directness. It yields a new topological proof,
which is constructive, transparent and economical.

Restricting our attention to the sub-operad of the operad of
cellular chains of normalized spineless cacti given by symmetric
top-dimensional cells  $CC_n^{top}(n)^{\mathbb {S}}$, we obtain a chain model for the operad
$\mathcal{G}Pl$ whose algebras are precisely graded pre-Lie algebras. Suitably shifting
degrees in this chain operad, we obtain the operad $\mathcal{P}l$ whose algebras
precisely  are pre-Lie algebras.

Let $L^*$ be a free $k$ (or $\mathbb{Z}$) module generated by an
element of degree $-1$.

\noindent {\bf Theorem \ref{preliecells}} {The operad
$CC_n^{top}(n)^{\mathbb {S}}\otimes k$ is isomorphic to the operad
$\mathcal{G}Pl$ for graded pre-Lie algebras. Furthermore the shifted
operad $(CC_n^{top}\otimes (L^*)^{\otimes E_w})^{\mathbb
{S}}(n)\otimes k$ is isomorphic to the operad $\mathcal{P}l$ for
pre-Lie algebras.

The analogous
statements also hold over $\mathbb{Z}$.}


The considerations above leading to the operations of the
Hochschild cochains are actually of a more general nature. To make
this claim precise, we analyze meta-structures on operads  and
show that their structure naturally leads to pre-Lie algebras,
graded pre-Lie algebras, Lie algebras and Hopf algebras. Here the
same ``flow-chart argument'' gives these structures a cell
interpretation. For another approach to relations between Hopf
algebras and co-operads, see \cite{LM,L}.

Specializing to the Hopf structure of the operad of the shifted
top-dimensional symmetric cells above and taking
$\mathbb{S}_n$-coinvariants, we obtain the renormalization
Hopf-algebra of Connes and Kreimer.

\noindent {\bf Proposition \ref{ckcor}} {\it $H_{CK}$ is isomorphic
to the Hopf algebra of $\Sn$ coinvariants of the sub-operad of
top-dimensional symmetric combinations of shifted cells
$((\CCcact)^{top})^{\mathbb{S}}\otimes (L^*)^{\otimes E_w}$ of the
shifted cellular chain operad of normalized spineless cacti
$CCcact\otimes (L^*)^{\otimes E_w}$.}


Our analysis thus unites  the pre-Lie definition of the
Gerstenhaber bracket in string topology, the arc operad and
the original work of Gerstenhaber with the renormalization
procedures of Connes and Kreimer in terms of natural
operations on operads.

Going beyond the algebraic properties of operads,  we prove that
any operad with a multiplication is an algebra over the cell model
of the little discs operad given by the cellular chains of
normalized spineless cacti. This is a generalization of Deligne's
conjecture to the operad level and realizes the  program of \cite{GV2}
by giving a moduli space interpretation to the brace algebra structure on an operad.

\noindent {\bf Theorem \ref{gendel}} {\it The generalized Deligne
conjecture holds. I.e.\  the direct sum of any operad algebra which
admits a direct sum  is an algebra over the chains of the little
discs operad in the sense that it is an algebra over the dg-operad
$\CCcact$}.

\vskip .5\baselineskip

The resemblance of the geometric realization of the chains
defining the homotopy Gerstenhaber structure and the algebraic
calculations of Gerstenhaber \cite{G} is striking, making the case
that cacti are the most natural topological incarnation of these
operations. In fact, using the translation formalism developed  in
this paper, the topological homotopies listed in \cite{KLP} are
the exact geometrization of the algebraic homotopies of \cite{G}.
As an upshot, the natural boundary map present in the spineless
cacti/arc description shows how the associative multiplication is
related to the bracket as a degeneration. In an algebraic
topological formulation this establishes that the pre-Lie
multiplication is basically a $\cup_1$ operation.

It is thus tempting to say that { the $\Arc$ operad  \cite{KLP}
is an underlying ``string mechanism'' for all of the above
structures.}

\subsection*{The paper is organized as follows:}

In the first paragraph, we introduce the types of trees we wish to
consider and several natural morphisms between them. This is
needed to fix our notation and allows us to compare our results
with the literature.

In the second paragraph, we recall the definitions of \cite{cact}
of the different types of cacti. Furthermore, we recall their
quasi-operad structure from \cite{cact} and the description of
spineless cacti as a semi-direct product of the normalized
spineless cacti and a contractible so-called scaling operad.
Finally, we recall the equivalence of the operad of spineless cacti  with
the operad of little discs.

In paragraph three, we give a cell
decomposition for the space of normalized spineless cacti.  We
also show that the induced quasi-operad structure on the chain
level is in fact an operad structure on the cellular chains of this
decomposition. This yields our chain model for the little discs
operad.

Paragraph four contains our new solution to Deligne's conjecture.
After recalling the definition of the Hochschild complex and the
brace operations, we  provide two points of view of the operation
of our cellular operad on the Hochschild cochains. One which is
close to the operation of the chains of the $\Arc$ operad on
itself and also to string topology and a second one which is based
on a description in terms of flow charts for substitutions and
multiplications among Hochschild cochains.

In the fifth paragraph, we show that the symmetric combinations of
these cells yield an operad which is the operad whose algebras are
precisely
graded pre-Lie algebras. We additionally show that the pre-Lie operad
has a natural chain interpretation in terms of the symmetric
combinations of the top-dimensional cells of our cell
decomposition for normalized spineless cacti. Furthermore, we
analyze the situation in which an operad that admits a direct sum
also has an element which acts as an associative multiplication.
In this setting, we prove a natural generalization of Deligne's
conjecture which states that a chain model of the little discs
operad acts on such an operad.

In the sixth paragraph, we use our previous analysis to define
pre-Lie and Hopf algebras for operads of $\mathbb{Z}$-modules or
any operad leading to $\mathbb{Z}$-modules. Applying the Hopf algebra
construction to  our chain
model operad for pre-Lie algebras
and then taking $\mathbb{S}_n$
co-invariants we obtain
the Hopf algebra of Connes and Kreimer.

In a final short paragraph, we comment on the generalization
to the $A_{\infty}$ case and the cyclic
case as well as on new developments.

\section*{Notation}

We denote by $\Sn$ the permutation group on n letters and by
$C_n$ the cyclic group of order n.

We denote the shuffles of two ordered finite sets $S$
and $T$ by $\Sh(S,T)$.  A shuffle of two finite ordered
sets $(S, \prec_S)$ and $(T \prec_T)$ is an order $\prec$ on
$S\amalg T$ which respects both the order of $S$ and that of
$T$, i.e.\ for $t,t' \in T$: $t \prec t'$ is equivalent to
$t\prec_Tt'$ and for $s,s' \in S$: $s \prec s'$ is equivalent to
$s\prec_Ss'$. We also denote by $\Sh'(S,T)$ the subset of
$\Sh(S,T)$ in which the minimal element w.r.t.\ $\prec$ is the
minimal element of $S$. Finally we denote the trivial shuffle in
which $s\prec t$ for all $s\in S$ and $t \in T$ by
$\prec_S\amalg\prec_T$

For any element $s$ of an ordered finite set $(S,\prec)$ which is
not the minimal element, we denote the element which immediately
precedes $s$ by $\prec(s)$. We use the same notation for a finite
set with a cyclic order.

We also fix $k$ to be a field of arbitrary characteristic.

We will tacitly assume that everything is in the super setting,
that is $\mathbb{Z}/2\mathbb{Z}$ graded. For all formulas, unless
otherwise indicated, the standard Koszul rules of sign apply.

We let $\Set,\Top,\Chain, \Vect$ be the monoidal categories
of sets, topological spaces, free Abelian groups and (complexes of) vector spaces over $k$
and call operads in these categories combinatorial, topological, chain
and linear operads, respectively.

\section{Trees}
In the following trees will play a key role for indexing purposes
and in the definition of operads and operadic actions.

\subsection{General Definitions}
\begin{df}
A graph $\G$ is a collection $(V(\G),F(\G),\delta:F(\G)\rightarrow V(\G), \imath:
F(\G)\rightarrow F(\G))$ with $\imath^2=id$ and no fixed points $\forall f\in F(\G):\imath(f)\neq f$.
 The set $V(\G)$ is called the set of vertices and
the set $F(\G)$ is called the set of flags. We let $E(\G)$ be the
set of orbits of $\imath$ and call it the edges of $\Gamma$.
Notice that $\delta$ induces a map $\del:E(\G)\rightarrow
V(\G)\times V(\G)$, and that the data $(V(\G),E(\G), \del)$
defines a CW-complex by taking $V(\G)$ to be the vertices or
0-cells and $E(\G)$ to be the 1-cells and using $\del$ as the
attaching maps. The realization of a graph is the realization of
this CW-complex.

A tree is a
graph whose realization is contractible.

A rooted tree is a tree with a marked vertex.

We call a rooted tree planted if the root vertex lies on a unique
edge. In this case we call this unique edge the ``root edge''.
\end{df}

We usually depict the root of a planted tree $\t$ by a small
square, denote the  root vertex by $root(\t)\in V(\t)$ and the
root edge by $e_{root(\t)}\in E(\t)$.

Notice that an edge $e$ of a graph or a tree  gives rise to a set
of vertices $\del (e)=\{v_1,v_2\}$. In a tree the set $\del e=
\{v_1,v_2\}$ uniquely determines the edge $e$. An orientation of
an edge is a choice of order of the two flags in the orbit. An
oriented edge is an edge together with an orientation of that
edge. On a tree giving an orientation to the edge $e$ defined by
the boundary vertices $\{v_1,v_2\}$  is equivalent to specifying
the either ordered set $(v_1,v_2)$ or the ordered set $(v_2,v_1)$.
If we are dealing with trees, we will denote the edge
corresponding to $\{v_1,v_2\}$ just by $\{v_1,v_2\}$ and likewise
the ordered edge corresponding to $(v_1,v_2)$ just by $(v_1,v_2)$.

An edge that has $v$ as a vertex is called an adjacent edge to
$v$. We denote by $E(v)$ the set of edges adjacent to $v$.
Likewise we call a flag $f$ adjacent to $v$ if $\delta(f)=v$ and
denote by $F(v)$ the set of flags adjacent to $v$. We call an
oriented edge $(f,\imath(f))$ incoming to $v$ if
$\delta(\imath(f))=v$. If $\delta(f)=v$ we call it outgoing.

Of course the oriented edges  are in 1-1 correspondence with the
flags by identifying $(f,\imath(f))$ with $f$, but it will be convenient to keep
both notions.

An edge path on a graph $\Gamma$ is an alternating sequence of
vertices and edges $v_1,e_1,v_2,e_2,v_3, \dots $ with $v_i\in
V(\Gamma), e_i\in E(\Gamma)$, s.t.\  $\del(e_i)= \{v_i,v_{i+1}\}$.
Notice that since we define this notion for a general graph, we
need to keep track of the vertices and edges.

\begin{df}
\label{contract}
 Given a tree $\t$ and an edge $e \in E(\t)$ one obtains a new
 tree by contracting the edge $e$. We denote this tree by
$\t/e$.

More formally let $e=\{v_1,v_2\}$, and consider the equivalence
relation $\sim$ on the set of vertices which is given by $\forall
w \in V(\t):w\sim w$ and $v_1 \sim v_2$. Then $\t/e$
 is the tree whose vertices are $V(\t)/\sim$ and whose edges are
 $E(\t) \setminus \{e\}/\sim'$ where $\sim'$ denotes
the induced equivalence relation $\{w_1,w_2\}\sim'\{w_1',w_2'\}$
if $w_1 \sim w_1'$ and $w_2 \sim w'_2$ or $w_2 \sim w_1'$ and
$w_1 \sim w'_2$.
\end{df}

\subsubsection{Structures on rooted trees}
A rooted tree has a natural orientation, toward the root. In fact,
for each vertex there is a unique shortest edge path to the root
and thus for a rooted tree $\tau$ with root vertex $root \in
V(\t)$ we can define the function $N:V(\t)\setminus
\{root\}\rightarrow V(\t)$ by the rule that
$$
N(v)=\text{the next vertex on unique path to the root starting at
$v$}
$$
This gives each edge $\{v_1,v_2\}$ with $v_2=N(v_1)$ the
orientation $(v_1,N(v_1))$.

We call the set $\{(w,v)|w \in N^{-1}(v)\}$ the set of incoming
edges of $v$ and denote it by $\In(v)$ and call the edge
$(v,N(v))$ the outgoing edge of $v$.

\begin{df} We define the arity of $v$ to be $|v|:=|N^{-1}(v)|$.
The set of leaves $V_{ leaf}$ of a tree is defined to be the
set of vertices which have arity zero, i.e.\ a vertex is a leaf if
the number of incoming edges is zero. We also call the outgoing
edges of the leaves the leaf edges and denote the collection of
all leaf edges by $E_{ leaf}$.
\end{df}

\vskip 1mm \noindent {\sc Caveat:} For $v\neq root$: $|v|=|E(v)|-1$. That is,
$|v|$ is the number of
incoming edges, which is the number of adjacent edges minus one.
For the root $|v|$ is indeed the
number of adjacent edges. \vskip 1mm

\begin{rmk}
For a rooted tree there is also a bijection which we denote by:
$\verttoedge: V(\t)\setminus \{root\} \ra E(\t)$. It associates to
each vertex except the root its unique outgoing edge $v\mapsto
(v,N(v))$.
\end{rmk}

\begin{df}
\label{branch} An edge $e'$ is said to be above $e$ if $e$ lies on
the edge path to the root starting at the vertex of $e'$ which is
farther from the root. The branch corresponding to an edge $e$ is
the subtree consisting of
  all edges which lie above $e$ (this includes $e$)
and their vertices. We denote this tree by $br(e)$.
\end{df}

\subsection{Planar trees}
\begin{df}
A planar tree is a pair $(\t,p)$ of a tree $\t$  together with a
so-called pinning $p$ which is a cyclic ordering of each of the
sets $E(v),v\in V(\t)$.
\end{df}

\subsubsection{Structures on planar trees}
A planar tree can be embedded in the plane in such a way that the
induced cyclic order from the natural orientation of the plane and
the cyclic order of the pinning coincide.

The set of all pinnings of a fixed tree is finite and is a
principal homogeneous set for the group
$$
\Stau := \times_{v \in V(\t)} \mathbb {S}_{|v|}
$$
where each factor  $\mathbb {S}_{|v|}$ acts by permutations on the
set of cyclic orders of the set $E(v)$. This action is given by an
identification of the symmetric group ${\mathbb S}_{|v+1|}$ with
the permutations of the set of the $|v|+1$ edges of $v$ and then
modding out by the subgroup of cyclic permutations which act
trivially on the set of induced cyclic orders of $E(v)$:
$\mathbb{S}_{|v|}\iso \mathbb {S}_{|v|+1}/C_{|v|+1}$.

\subsubsection{Planted planar trees}
\label{pptrees}
Given a rooted planar tree $(\t,root)$ there is a linear order on
each of the sets $E(v)$, $v\in V(\t)\setminus \{root\}$. This
order is given by the cyclic order and designating the outgoing
edge to be the smallest element. The root vertex has only a cyclic
order, though.

Since the root of a planar planted tree  has only one incoming
edge and no outgoing one such a tree has a linear order at all of
the vertices. Vice-versa providing a linear order of the edges of
a root vertex is tantamount to planting a rooted tree by adding a
new root edge which induces the given linear order. Therefore
these are equivalent pictures and we will use either one of them
depending on the given situation.

Furthermore, on such a tree there is an edge  path which passes through
all the edges exactly twice ---once in each direction--- by starting
at the root going along the root edge and at each vertex
continuing on the next edge in the cyclic order until finally
terminating in the root vertex. We call this path the outside
path.

By omitting recurring elements, that is counting each vertex or
edge only the first time it appears, the outside path endows the
set $V(\t)\amalg E(\t)$ with a linear order $\prec^{(\t,p)}$. The
smallest element is the root edge and the largest element in this
order is the root vertex. This order induces a linear order on the
subset of vertices $V(\t)$, on the subset of all edges $E(\t)$, as
well as a linear order on each of the subsets $E(v)$ for all the
vertices. In this order on $E(v)$ the smallest element is the
outgoing edge. This order will be denoted by $\prec_v^{(\t,p)}$.
We omit the superscript for $\prec^{(\t,p)}$ if it is clear from
the context.

\subsubsection{Labelled trees}

\begin{df}
For a finite set $S$ an $S$-labelling for a tree is an injective
map $\lb: S\rightarrow V(\t)$. An $S$-labelling of a tree yields a
decomposition into disjoint subsets of $V(\t)= V_l \amalg V_u$
with $V_l=\lb(S)$. For a planted rooted tree, we demand that the
root is not labelled: $root \in V_u$.

An $n$-labelled tree is a tree labelled by $\bar
n:=\{1,\dots,n\}$. For such a tree we call $v_i:=\lb(i)$.

A fully labelled tree $\t$ is a tree such that $V_l=V(\t)$.
\end{df}

\subsection{Black and white trees}
\begin{df}
A black and white graph  (b/w graph) $\G$ is a graph together with
a function $\color:V(\G) \rightarrow \{0,1\}$.

We call the set $V_w(\t):=\color^{-1}(1)$ the set of white
vertices and call the set $V_b(\t):=\color^{-1}(0)$ the set of
black vertices.

By a bipartite b/w tree we understand a b/w tree whose edges only
connect vertices of different colors.

A $S$-labelled b/w tree is a b/w tree in which exactly the white
vertices are labelled, i.e.\  $V_l=V_w$ and $V_u=V_b$.

For a rooted tree we call the set of black leaves the tails.

A rooted b/w tree is said to be without tails if all the leaves
are white.

A rooted b/w tree is said to be stable if there are no black
vertices of arity 1, except possibly the root.

A rooted b/w tree is said to be fully labelled if all vertices
except for the root and the tails are white and labelled.

In a planar planted rooted b/w tree, we require that the root
be black.
\end{df}

\begin{df}
For a black and white bipartite tree, we define the set of white
edges $E_w(\t)$ to be the edges $\{b,N(b)\}$ with $N(b)\in V_w$
and call the elements white edges. Likewise we define
$E_b=\{\{w,N(w)\}|N(w)\in V_b\}$ with elements called black edges,
so that there is a partition $E(\t) = E_w(\t) \amalg E_b(\t)$.
White edges thus point towards white vertices and black edges
towards black vertices in the natural orientation towards the
root.
\end{df}

\begin{nota}
For a  planar planted b/w tree, we understand the adjective
bipartite to signify the following attributes:
\begin{enumerate}
\item both of the
vertices of the root edge are black, i.e.\ $root$ and the vertex
$N^{-1}(root)$ are black
\item the tree after deleting the tail vertices
and their edges, but
keeping the other vertices of these edges, is bipartite otherwise.
\end{enumerate}

The root edge is
considered to be a black edge. Also in the presence of tails, all
non-white tail edges are considered to be black.
\end{nota}

\begin{df} For a planar planted b/w bipartite tree  $\t$ and a white  edge
$e=(b, N(b))$ we define the branch of $e$ denoted by $br(e)$  to
be the planar planted bipartite rooted tree  given as follows.
\begin{enumerate}
\item
The vertices and edges are those of the branch of $e$ as defined
in Definition \ref{branch}.
\item
The colors of all vertices except
$N(b)$ are kept and the color of $N(b)$ is changed to black. This
black vertex is defined to be the root.
\item In the case that the tree
$\t$ is also labelled, $br(e)$ is considered to be labelled by the
set of labels of its white vertices --- we stress that this does
not include the root $N(b)$, which is  the root of $br(e)$.
\end{enumerate}
\end{df}

\subsubsection{Notation I} N.B.\ A tree can have several of the attributes mentioned above;
for instance, we will look at  bipartite planar planted rooted
trees. To fix the set of trees, we will consider the following
notation. We denote by $\mathcal T$ the set of all trees and use
sub- and superscripts to indicate the restrictions. The
superscript $r,pp,nt$ will mean rooted and planar planted, without
tails while the subscripts $b/w,bp,st$ will mean black and white,
bipartite, and stable, where bipartite and stable insinuate that
the tree is also b/w. E.g.\

\begin{tabular}{ll}
 $\rtree$&The set of all rooted trees\\
 $\bwpptree$&The set of  planar planted b/w trees\\
 $\bppptree$&The set of planar planted  bipartite trees\\
 $\stpptree$& The set of planar planted stable b/w trees\\
\end{tabular}

Furthermore we use the superscripts $fl$
for fully
labelled
trees. E.g.\

\begin{tabular}{ll}
 $\flrtree$&The set of all rooted fully labelled trees\\
 \end{tabular}

 We furthermore use the notation
that $\mathcal{T}(n)$ denotes the $n$-labelled trees and adding
the sub and superscripts denotes the $n$-labelled trees of that
particular type conforming with the restrictions above for the
labelling. Likewise $\mathcal{T}(S)$ for a set $S$ are the $S$-labelled
trees conforming with the restrictions above for the
labelling. E.g.\

\begin{tabular}{ll}
 $\bwpptree(n)$&The set of planar planted b/w trees with $n$  white vertices \\
 &which are labelled by the set $\{1,\dots,n\}$.\\
\end{tabular}

\subsubsection{Notation II}
Often we wish to look at the free Abelian groups or free vector
spaces generated by the sets of trees. We could introduce the
notation $Free(\mathcal{T}, \Z)$ and $Free(\mathcal{T}, k)$ with
suitable super- and subscripts, for the free Abelian groups or
vector spaces generated by the appropriate trees. In the case that
there is no risk of confusion, we will just denote these freely
generated objects again by  $\mathcal{T}$ with suitable sub- and
superscripts  to avoid cluttered notation. If we define a map on
the level of trees it induces a map on the level of free Abelian
groups and also on the $k$-vector spaces. Likewise by tensoring
with $k$ a map on the level of free Abelian groups induces a map
on the level of vector spaces. Again, we will denote these
maps in the same way.

\subsubsection{Notation III}
If we will be dealing with operads of trees, we will consider the
collection of the $\mathcal{T}(n)$ with the appropriate sub- and
superscripts. Again to avoid cluttered notation when dealing with
operads, we  also  denote the whole collection of the
$\mathcal{T}(n)$ just by $\mathcal{T}$ with the appropriate sub-
and superscripts.

\begin{df} There are some standard trees, which are essential in our
study, these are the n-tail tree $l_n$, the white $n$-leaf tree
$\t_n$, and the black $n$-leaf tree $\t^b_n$, as shown in Figure
\ref{nleaf}.
\begin{figure}
\epsfbox{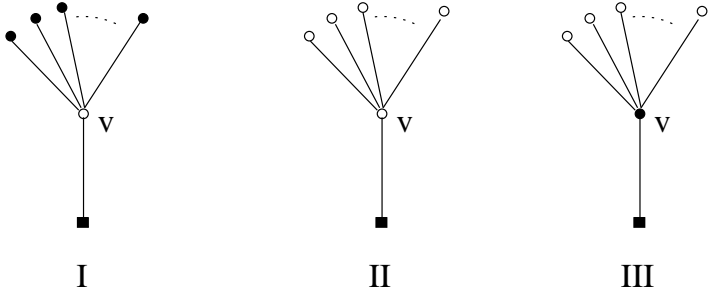} \caption{\label{nleaf}I. The n-tail tree
$l_n$, II. The  white $n$-leaf tree $\t_n$, III. The black
$n$-leaf tree $\t^b_n$.}
\end{figure}
\end{df}

\subsection{Maps between different types of trees}

\subsubsection{The map $\cppin: \rtree \rightarrow \wlbptree$}

First notice that there is a map from planted trees to rooted trees given by
contracting the root edge. This map actually is a bijection
between planted and rooted trees. The inverse map is given by
adding one additional vertex which is designated to be the new
root and introducing an edge from the new root to the old root. We
call this map $\plant$.

Secondly, there is a natural map $\pin$ from the free Abelian group of planted
trees to that of planted planar trees. It is given by:
$$
pin(\t) = \sum_{\text {p $\in$ Pinnings}(\t)} (\t,p)
$$

Finally there is a map from planted planar trees to planted planar
bipartite trees without tails. We call this map $\bp$. It is given
as follows. First color all vertices white except for the root
vertex which is colored black, then insert a black vertex into
every edge.

In total we obtain a map
$$
\cppin:= \bp \circ \pin \circ \plant: \rtree  \rightarrow
\wlbptree
$$
that plants, pins and colors and expands the tree in a bipartite
way.

Using the map $\cppin$, we will view $\rtree$ as a subgroup of
$\wlbptree$. The image of $\rtree$ coincides with the set of
invariants of the actions $\mathbb{S}(\t)$. We will call such an
invariant combination a symmetric tree.

\begin{rmk}
The inclusion above extends to an inclusion of the free Abelian
group of fully labelled rooted trees to labelled bipartite planted
planar
 trees: $\cppin:\flrtree(n) \rightarrow \wlbptree(n)$.
\end{rmk}

\subsubsection{The map $\forgetainfty :\stpptree \rightarrow \bppptree$}

We define a map from the free groups of stable b/w planted planar
trees to the free group of bipartite b/w planted planar trees in
the following way:
 First, we set to zero any tree which has black vertices
whose arity is greater than two. Then, we contract all edges
which join two black vertices. And finally, we insert a black
vertex into each edge joining two white vertices. We call this map
$\forgetainfty$.

Notice that $\forgetainfty$ preserves the condition of having no
tails and induces a map on the level of labelled trees.

This nomenclature is  chosen since this map in a
sense  precisely forgets the trivial $A_{\infty}$ structure of an associative
algebra in which all higher multiplications are zero and
all $n$-fold iterations of the multiplication agree.

\subsection{An operad structure on $\bwbptree$}

\subsubsection{Grafting planar planted b/w trees at leaves}
\label{leafgraftingpar}
Given two trees, $\t \in \bppptree(m)$, $\t' \in \bppptree(n)$ and
a white vertex $v_i$ which is a leaf of $\t$, we define $\t
\circ_i \t'$ by the following procedure:

First identify the root of $\t'$ with the vertex $v_i$. The image
of $v_i$ and $root(v_i)$ is taken to be black and unlabelled. The
linear order of all of the edges is given by first enumerating the
edges of $\t$ in their order until the outgoing edge of $v_i$ is reached,
then enumerating the edges of
the tree $\t'$ in their order and finally the rest of the edges of $\t$ in their
order, the latter being all the edges following the outgoing edge of $v_i$ in the order
of $\t$.

Second contract the image of the root edge of $\t'$, i.e.\ the
image of the edge $e_{root}(\t')$ under the identification $v_i
\sim root(\t')$, and also contract the image of the outgoing edge of $v_i$, i.e.
the image after gluing and contraction of the edge $(v_i,N(v_i))$.

The root of this tree is specified to be the image of the root of
$\t$ and the labelling is defined in the usual operadic way. The
labels $1,\dots, i-1$ of $\t$ are unchanged, the labels
$1,\dots,n$ of $\t'$ are changed to  $i, \dots, n+i-1$, and finally the
labels $i+1, \dots , m$ of $\t$ are changed to  $i+n, \dots, m+n-1$.

\subsubsection{Cutting branches}
Given a tree $\t\in \bppptree(n)$ and a vertex $v_i$ we denote the
tree obtained by cutting off all branches at $v$ by $ \cut(\t,v).
$ This is the labelled subtree of $\t$  consisting of all edges
below (viz.\ not above) the incoming edges  of $v$ and the
vertices belonging to these edges. We stress that the outgoing
edge $\verttoedge(v)$ is a part of this tree as is $v$. By keeping
the labels of the remaining white vertices the tree $\cut(\t(v)$
becomes an $S$-labelled planar planted bipartite tree. Here
$S\subset \{1,\dots,n\}$ is the set of labels of the subtree under
consideration.

Let $\br(\t,v)$ denote the set of the branches of the
incoming edges of $v$ which is ordered by the order $\prec$ induced by  $\prec^{\t}_v$.
I.e.\ $br(e_i) \prec br(e_j)$ if and only if $e_i\prec_v^{\t} e_j$.

$$
\br(\t,v)=(\{br(e_i)|e_i \in \In(v_i)\},\prec)
$$

\subsubsection{Grafting branches}

Fix a tree $\t \in \bppptree(S_0)$ and an ordered set of trees
$\t'_i \in \bppptree(S_i): i\in\{ 1,\dots, m\}$. Let
$R:=(\{e_{root}(\t'_1), \dots, e_{root}(\t'_m)\}, \prec_R)$ be the
ordered set in which $e_{root}(\t_i) \prec_R e_{root}(\t_j)$ if
and only if $i<j$.

Set $E:=(E(\t)\setminus \{e_{root}\},\prec^{\t})$ and define
$\whitevert:E\rightarrow V(\t)$ to be the map which maps each edge
in $E$ to its unique white vertex. Denote the minimal element of
$E$ by $e_{min}$. This is the edge which immediately follows
$e_{root}(\t)$ in the linear order of $\t$.

Given a shuffle $\prec \in \Sh'(E,R)$ such that the minimal
element of the ordered set $(E\amalg R,\prec)$ is $e_{min}$ we
define the grafting of the branches $\t'_1,\dots,\t'_m$ onto $\t$
with respect to $\prec$ to be the labelled b/w tree obtained as follows:

\begin{enumerate}
\item
Identify the roots of the $\t'_i$ with the white vertex of
the edge immediately preceding the root edge: $root_{\t'_i}\sim
w=\whitevert(\prec(e_{root(\t')})$

\item Designate the image to be
a white vertex with the label $L^{-1}(w)$.

\item Endow this tree with the planted planar structure induced by
the order $\prec$ together with the orders $\prec^{\t'_i}$ and
$\prec^{\t}$. The root of  this tree is the image of the root of
$\t$.
\end{enumerate}
This tree is again in $\bppptree$ and is labelled by $S_0\amalg S_1\amalg
\dots \amalg S_M$. We denote it by

$$
\gr(\t;\t_1,\dots,\t_n;\prec)
$$

\subsubsection{Signs}
For a shuffle of sets with weighted elements $wt_S:S\rightarrow
\mathbb{N}$ $wt_T:T\rightarrow \mathbb{N}$, we define the sign of
the shuffle to be the sign obtained from shuffling the elements
of $T$ past the elements of $S$, i.e.\

$$\sign(\prec)= \prod_{t\in T}(-1)^{\sum_{s\in S: t\prec s}wt_S(s)wt_T(t)}.$$

We define the weight function  by
\begin{equation}
\label{wtdef}
wt(e) = \begin{cases} 1& \text{if $e$ is white}\\
0& \text{if $e$ is black and $e$ is not a root edge}\\
|E_w(\t)|&\text{if $e$ is the root edge of $\t$}
\end{cases}
\end{equation}

\begin{df}
We define
 the grafting of the $\t_i$
onto $\t$ as branches to be the signed sum over all possible graftings
using the weight function
(\ref{wtdef}):

\begin{equation}
\label{treegraft}
 \gr(\t;\t_1,\dots,\t_n):= \sum_{\prec \in
\Sh'}\sign(\prec) \gr(\t;\t_1,\dots,\t_n;\prec)
\end{equation}
\end{df}

\subsubsection{An operad structure for $\bppptree$}
\label{insertionexamples}

With the above procedures, we define operadic compositions $\t\circ_i\t'$  for
$\bppptree$ as follows: first cut off the branches corresponding
to the incoming edges of $v_i$. Second graft $\t'$ as a planar
planted tree onto the remainder of $\t$ at the vertex $v_i$ which
is now a leaf according to the grafting procedure of \S\ref{leafgraftingpar}.
Finally sum over the possibilities to graft the cut
off branches onto the white vertices of the resulting planar tree
which before the grafting belonged to $\t'$. Here one only sums
over those choices in which the order of the branches given by the
linear order at $v_i$ is respected by the grafting procedure.
I.e.\ the branches after grafting appear in the same order on the
grafted tree as they did in $\t$.

\begin{figure}
\epsfxsize = \textwidth \epsfbox{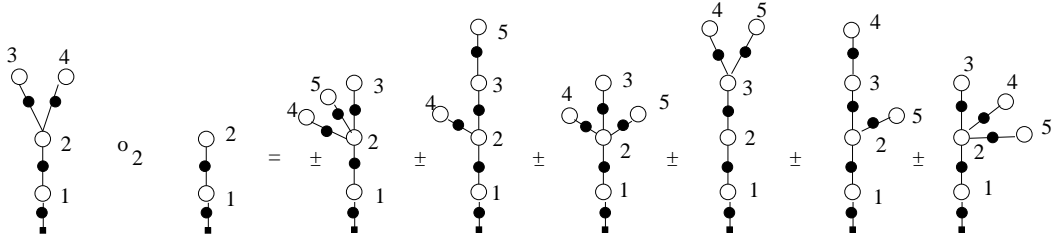}
\caption{\label{bptreeins} Example of the insertion of a bipartite
planted planar tree}
\end{figure}

An example of such an insertion is depicted in figure
\ref{bptreeins}.

\begin{df}
Given $\t \in \bppptree(m)$ and  $\t'\in \bppptree(n)$, we define
the tree $\t\circ_i \t' \in \bppptree(m+n-1)$
\begin{equation}
\label{treecircs}
 \t \circ_i \t' := \gr(\cut(\t,v_i)\circ_i
\t';\br(\t,v_i))
\end{equation}
with the following relabelling: the labels $1,\dots, i-1$ of
vertices which formerly belonged to $\t$ remain unchanged. The
labels of the vertices $1,\dots, n$ of the vertices which formerly
belonged to $\t'$ are relabelled $i, \dots, i+n-1$ and the
remaining vertices of those which formerly belonged to $\t$ which
used to be labelled by $i+1,\dots,m$ are now relabelled by
$i+n,\dots,m+n-1$.
\end{df}

\subsubsection{Labelling by sets}
\label{setlabelleing}
 There is a way to avoid labelling and
explicit signs by working with tensors and operads labelled by
arbitrary sets (\cite{D,MSS,KS}). In this case, if $S$ and $S'$
are the indexing sets for $\t$ and $\t'$ and $i\in S$, then the
indexing set of $\t \circ_i \t'$ is given by $S \setminus \{i\}
\amalg S'$. To obtain the signs, one associates a free
$\mathbb{Z}$-module (or $k$ vector space) generated by an element
of degree minus one to each white edge. See Definition
\ref{shifted}.

\subsubsection{Positive signs}
We define $\t \circ_i^+ \t'$ just as above only with the
$wt(e)\equiv 0$, that is as the formal sum with the same summands
and only positive coefficients.

\subsubsection{Contracting trees}
Given a rooted subtree $\t' \subset \t \in \wlbptree$ with white
leaves and black root, we define $\t/\t'$ to be the tree obtained
by collapsing the subtree $\t'$ to one black edge by identifying
all white and all black vertices of $\t'$. I.e.\ let
$v\sim_V^{\t'}v'$ if $v=v'$ or $v,v'\in V(\t')$ and
$\color(v)=\color(v')$
 then $V(\t/\t')/\sim_V^{\t'}$. Likewise let $e
\sim_E^{\t'} e'$ if $e,e'\in E(\t')$ or $e=e'$ and set
$E(\t/\t')=E(\t)/\sim_E^{\t}$.

For any tree labelled tree $\t\subset \wlbptree(n)$ with labelling
$\lb$, we set $\t^{+i}$ to be the same underlying tree, but with
the shifted labelling function $\lb^{+i}:\{i, \dots, i+n-1\}
\rightarrow E_w(\t)$ given by $\lb^{+i}(k)=\lb(k-i+1)$.

By abuse of notation we use $\t'$ to denote a subtree $\t' \subset
\t$  and  the planted planar tree obtained from $\t'$ by planting
the root vertex while preserving the linear order already
present on the original edges of $\t'$. Vice-versa given a planted tree $\t$,
when identifying it with a subtree the root edge will be contracted,
but the linear order of the subtree has to coincide with that of
the tree considered as ``free standing''.

Consider $\t'\in \wlbptree(m)$. If $\t^{\prime +i}\subset \t\in
\wlbptree(m+n-1)$ as a labelled subtree, we write $\t' \subset_i
\t$. In this case we label $\t/\t'$ as follows. The labels $1,
\dots, i-1$ remain unchanged, the label of the vertex representing
the white vertices of the contracted $\t'$ is set to be $i$ and
the labels $i+m,\dots ,n+m-1$ are changed to $i+1, \dots, n$. Set
\begin{equation}
T(\t,\t',i):=\{\tilde \t|\t' \subset_i \tilde \t \text{ and }
\tilde \t/\t'=\t\}
\end{equation}

\begin{rmk}
\label{refomurem} With the above definitions, we can rewrite
equation (\ref{treecircs}) as
\begin{equation}
\label{newtreecircs} \t \circ_i \t' := \sum_{\tilde \t \in
T(\t,\t',i)}
 \sign(\prec^{\tilde \t}) \tilde \t
\end{equation}
 Where $\prec^{\tilde \t}$ is considered as the shuffle
 $(E_w(\tilde \t)=E_w(\t)\amalg E_w(\t'),\prec)$ of the
ordered sets $(E_w(\t),\prec^{\t})$ and $(E_w(\t'),\prec^{\t'})$.
Notice that the compatibility of the orders is automatic.
\end{rmk}

\begin{prop}
\label{treeopprop}
 The gluing maps (\ref{treecircs}) (with or without signs) together with
the symmetric group actions permuting the labels turn
$\bppptree:=\{\bppptree(n)\}$ into an operad.
\end{prop}

\begin{proof} This is a straightforward calculation especially in view
of the reformulation of \S \ref{refomurem}. An alternative
topological way to prove the associativity of $\circ^+_i$ is given
in Corollary \ref{topasscor}. The signs for the maps $\circ_i$ then follow from
\S\ref{signchoices} using Definition \ref{shifted}.
\end{proof}

\begin{rmk}
With the above compositions $\wlbptree$ is a suboperad of
$\bwbptree$.
\end{rmk}

\subsubsection{The differential on $\wlbptree$}
\label{diffpar}
 There is a  differential on  $\wlbptree$ which we
will now define in combinatorial terms. We will show later that
 it has a natural interpretation as the differential
 of a cell complex.

Recall that for a planted planar tree there is a linear order on
all edges and therefore a linear order on all subsets of edges.

\begin{df} Let $\t\in \wlbptree$. We set $E_{ angle}=E(\t)\setminus
(E_{ leaf}(\t)\cup \{e_{root}\})$ and we denote by
$\num_E:E_{ angle} \rightarrow \{1,\dots,N\}$  the bijection
which is induced by the linear order $\prec^{(\t,p)}$.
\end{df}

\begin{figure}
\epsfxsize = \textwidth \epsfbox{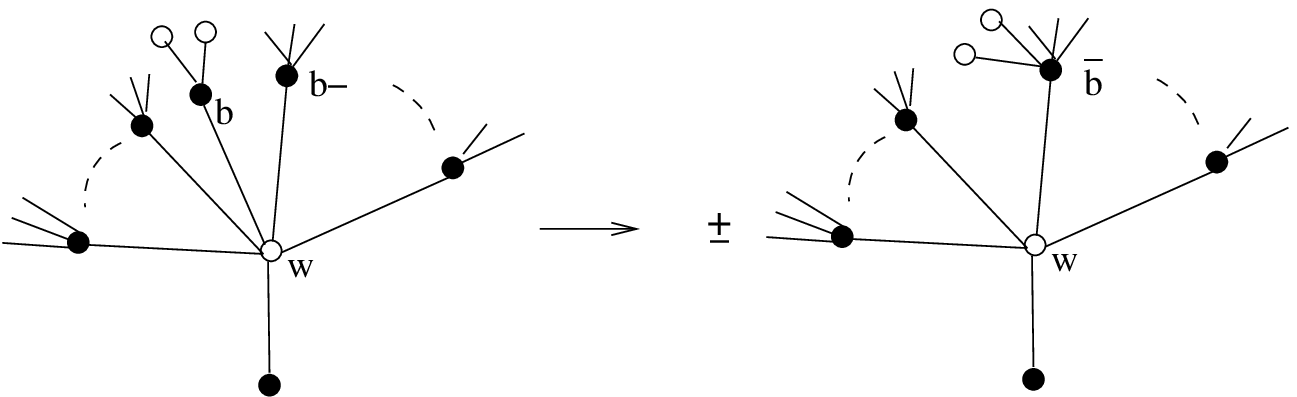} \caption{\label{del} The
tree $\del_e(\tau)$}
\end{figure}

\begin{df}
\label{diffdef}
 Let $\t\in \wlbptree$, $e\in E_{ angle}$,
$e=\{w,b\}$, with $w\in V_w$ and $b\in V_b$. Let $e-=\{w,b-\}$ be
the edge preceding $e$ in the cyclic order $\prec^{\t}_w$ at $w$.
Then $\del_e(\tau)$ is defined to be the planar tree obtained by
collapsing the angle between the edge $e$ and its predecessor in
the cyclic order of $w$ by identifying $b$ with $b-$ and $e$ with
$e-$.
 Formally

\begin{eqnarray*}
w=\whitevert(e), \quad e-=\prec^{\t}_w(e),&& \{b-\}= \del(e-)\cap V_b(\t)\\
V_{\del_e(\tau)}=V(\t)/(b\sim b-),&&
E_{\del_e(\tau)}=E_{\tau}/(e\sim e-)
\end{eqnarray*}
The linear order of $\del_e(\t)$ is given by keeping the linear
order at all vertices which are not in the image $\bar b$ of $b$
and $b-$ and using the linear order
$$(\In(\bar b), \prec_{\bar b}^{\del_e(\t)})
=(\In(b-)\amalg\In(b), \prec^{\t}_{b-}\amalg \prec^{\t}_{b})
$$
extended to $E(\bar b)$ by declaring the image of $e$ and $e-$ to
be the minimal element.
See Figure \ref{del} for an example.

\end{df}
\begin{df}
We define the operator $\del$ on the space $\wlbptree$ to be given
by the following formula
\begin{equation}
\label{treedifferential}
 \del(\t) := \sum_{e\in E_{ angle}}
(-1)^{\num_E(e)-1} \del_e (\tau)
\end{equation}
\end{df}

Denote by $\wlbptree(n)^k$ the elements of $\wlbptree(n)$ with $k$
white edges.

\begin{prop} The map $\del:\wlbptree(n)^k \rightarrow \wlbptree(n)^{k-1}$
is a differential for $\wlbptree$ and turns $\wlbptree$ into a
differential operad.
\end{prop}

\begin{proof}
The fact that $\del$ reduces the number of white edges by one is
clear. The fact that $\del^2=0$ follows from a straightforward
calculation. Collapsing two angels in one order contributes
negatively with respect to the other order. The compatibility of
the multiplications $\circ_i$ is also straightforward.
\end{proof}

\subsubsection{Other choices of signs}
\label{signchoices} The way the signs are fixed in the above
considerations is
 by giving the white edges the weight 1 and the black vertices
 the weight 0. Once the order of the edges for trees is fixed
 all signs are  the standard signs obtained
 from permuting weighted (i.e.\ graded) elements.
 We chose the ``natural'' order in which the edges are enumerated
 with respect to $\prec^{\t}$, i.e.\ the order derived from the
 embedding into the plane.

Another choice of ordering would be ``operadic'' in which the
white edges are enumerated first according to the label of their
incident white vertex and then according to their linear order at
that vertex. We leave it to the reader to make the necessary
adjustments in the formulae (\ref{treegraft}), (\ref{treecircs})
and (\ref{treedifferential}) to adapt the signs to this choice.

Finally one can avoid explicitly fixing an order if one works with
operads over arbitrary sets (see also \S \ref{setlabelleing}).

\subsubsection{Other tree insertion operads and compatibilities}
There are two tree insertion operads structures already present in
the literature on rooted trees $\rtree$ \cite{CL}, or to be more
precise on $\flrtree$ and, historically the first, on planar
planted stable b/w trees without tails $\wlptree$ \cite{KS}.

In the gluing for $\flrtree$ one simply omits mention of the order. And in
the case of $\wlptree$ one also allows gluing to the images of the
black vertices of $\t'$. Also in the  case of $\flrtree$ the basic
grafting of trees is used (no contractions), while in the case of
$\wlptree$ the grafting for planted trees is used, i.e.\ the image
of the root edge is contracted, but not the outgoing edge of
$v_i$.

The signs for the first gluing are all plus \cite{CL} and in the
second gluing they are given by associating the weight $1$ to all
edges except the root and weight $-2$ to the white vertices and
weight $2$ to the root vertex. The latter of course do not
contribute to the signs.

\begin{prop}
\label{symmetrictreesprop}  The map
$\forgetainfty:\wlptree\rightarrow \wlbptree$ is an operadic map.
Moreover,  it  is  a map of differential graded operads, if one
reverses the grading of $\wlbptree$, i.e.\ defining
$(\wlbptree)(n)^k$ to have degree $-k$.
\end{prop}

\begin{proof}
The operadic properties  of $\forgetainfty$ over $\Z/2\Z$ follow
from the fact that the summands given by gluing branches to black
vertices in the composition of \cite{KS} are sent to zero by the
map $\forgetainfty$. In the calculation for the differential the
unwanted terms disappear by the same argument. The compatibility
of signs follows from the fact that the trees that survive
$\forgetainfty$ have trivalent black subtrees. These have an odd
number of edges counting all the edges incident to the black
vertices. Let the outgoing edge be the outgoing edge of the lowest
vertex of the subtree. Collapsing the subtree then corresponds to
assigning an odd weight to the outgoing edge and weight zero to
the incoming edges, since there is no gluing to black vertices. In
the case of edges connecting two white vertices, one inserts a
black vertex and the weight can be transferred to the white edge,
again since there is no gluing onto black vertices. Now the parity
of the induced weights on the edges coincides with the weights  we
defined on $\wlbptree$. The weights of the vertices being even
play no role. Hence the signs agree. In fact taking the signs of
the vertices into account, we would exactly obtain weight $-1$ for
white edges and weight $0$ for black ones. The compatibility of
the degrees follows from the fact that in $\wlptree$ a black
vertex $b$ will contribute degree $-(|b|-2)$ and a white vertex
$w$ will contribute degree $-|w|$ to the total degree. Thus in the
case in which the tree has only binary black vertices the total
degree is the sum over the $-|w|$ which is the negative of the
grading in $\wlbptree$.
\end{proof}

\begin{df}
\label{shifted} We define the shifted complex $S^+\wlptree$ to be
given by $S^+\wlptree(n):=\wlptree(n)\otimes
 L^{\otimes E_w}$
where $L$ is a freely generated $k$ (or $\mathbb{Z}$ module)
generated by an element $l$ of degree $1$ and we
 used the notation of tensor products indexed by sets.
This means that generators are given by
 $\t\otimes l^{\otimes E_w(\t)}$.  The differential is the tree
differential without signs on the trees which just collapses the
angles and the multiplications on the component of trees are the
$\circ_i^+$. Now the signs come from the tensor factors and their
permutations as induced by the maps \S \ref{setlabelleing} and
$L^{\otimes E_w(\t)}\rightarrow L^{\otimes E_w(\del_e(\t))}$. Here
the latter map can be induced by the multiplication map
$\mu:L_e\otimes L_{e-}\rightarrow L_{\bar e}$ defined by
$\mu(l\otimes l)= l$.

It is clear that there is an operadic
isomorphism $S^+\wlbptree \simeq (\wlbptree,\circ_i)$.

Notice that if we shift back with the dual line $L^*$, i.e.\ set
$S\wlbptree=\wlbptree \otimes (L^*)^{\otimes E_w}$ where now the
compositions on the tree factor are given by $\circ_i$, then there
is an operadic isomorphism $S\wlbptree\simeq
(\wlbptree,\circ_i^+)$.

In an analogous way we define $\flrtree \otimes L^{\otimes E}$.
\end{df}

\begin{prop}
\label{cppinop} The map $\cppin:\flrtree\rightarrow \wlbptree$ is
injective and its image are the symmetric combinations of the
trees of $\wlbptree$ with $|E_w|=|V_w|-1$ which we denote by
$((\wlbptree)^{top})^{\mathbb{S}}$. Moreover
\begin{itemize}
\item[(a)] $\cppin:\flrtree\rightarrow\wlbptree$ is an operadic
embedding into the operad $(\wlbptree,\circ_i^+)$.  This
corresponds to assigning weight zero to all edges in the formalism
above.

\item[(b)] The  map $\cppin$ also induces an operadic  embedding
of $\flrtree\otimes L^{\otimes E}$
 into the  $(\wlptree,\circ_i)$, via
$$
 \begin{CD}
 \flrtree \otimes
 L^{\otimes E}
 @>\cppin \otimes id^{\otimes |E(\t)|}>>
 \wlbptree \otimes L^{\otimes E_w}=S^+\wlbptree
 \end{CD}
$$
We will also denote this map by $\cppin$.
\end{itemize}
\end{prop}
\begin{proof}
For rooted trees the composition defined in \cite{CL} tells us to
cut off the branches, glue the second tree to the truncated tree
and redistribute the branches. Now the map $\cppin$ associates to
each tree a sum whose summands are uniquely determined by a linear
order on the underlying rooted tree, which is obtained by
forgetting the linear order, contracting the root edge and all
black edges. The possible linear orders on the composed tree are
naturally in a 1-1 correspondence between the linear orders on the
truncated tree, the tree which is grafted on and a compatible
order of the re-grafted branches. These are the terms appearing in
the gluing of the planar planted trees. Before and after the
embedding, the symmetric group actions produces no signs. The
second statement follows clearly from the first.
\end{proof}

\section{Species of Cacti and their relations to other operads}
\label{cactchapter}

\subsection{Spineless Cacti}
In this section, we review the spaces of spineless cacti,
normalized spineless cacti, their relation to each other and the
little discs operad. These spaces were introduced in
\cite{KLP,cact}, but in order to facilitate the reading we
reiterate one definition of the sets corresponding to spineless
cacti and normalized spineless cacti and define operadic
respectively quasi-operad maps on these sets.  In order to not
disrupt the flow of the paper too much, we relegate some of the
technical details to the Appendix. To be completely
self-contained, we also define a topology on these sets which
makes spineless cacti into an operad of spaces and normalized
spineless cacti into a homotopy associative quasi-operad. After
this we briefly recall other equivalent ways of giving a topology.

\subsubsection{Background on spineless cacti and cacti}
In \cite{cact} we introduced the operad of spineless cacti. It is
a suboperad of the operad of cacti which was introduced
 by Voronov in \cite{Vor} descriptively as tree-like
configurations of circles in the plane. One precise definition of
the topology of cacti was given in \cite{KLP} via an operadic
embedding of cacti into the operad $\DArc$ which is
a deprojectivised version of the $\Arc$ operad \cite{KLP} .
 By definition $\DArc=\Arc \times \mathbb{R}_{>0}$. It carries
 the obvious  action of $\mathbb{R}_{>0}$ acting freely on the right
factor. This action is a global re-scaling. The spaces $\Arc(n)$
of the $\Arc$ operad are  open subsets of a cell complex
\cite{KLP}. This realizes cacti as a subset of a cell complex
crossed with $\mathbb{R}_{>0}$ whence it inherits the subspace
topology. Equivalently, using the map $\Loop$ of \cite{KLP} one
arrives at a description of the set of cacti as marked tree-like
ribbon graphs with a metric. Again, one obtains a subspace
topology from that of marked metric ribbon graphs. A marking is
the choice of a point on each cycle. In this setting tree-like
means that the ribbon graph is of genus zero and that there is one
marked cycle which passes through all of the un-oriented edges. A
very brief summary is given below. A short but detailed summary of
both these constructions is given in Appendix B of \cite{cact}.

\subsubsection{Background on the normalized versions}
In \cite{cact} we also introduced  normalized spineless cacti and
normalized cacti. The condition of normalization means in the
description via configurations of circles that all the circles are
of length one. The normalized versions are homotopy equivalent to
the non-normalized versions as spaces. They are not operads, but
they can be endowed with a quasi-operad structure which is
quasi-isomorphic to the operad structure of the non-normalized
version, i.e.\ it agrees with the operad structure of cacti
respectively spineless cacti on the level of homology \cite{cact}.

The main result of the next section is that there is a cell
decomposition for normalized spineless cacti, such that the
induced quasi-operad structure on the level of cellular chains is
an operad structure, and hence gives an operadic chain model for
spineless cacti.

We start by recalling the definition of spineless cacti and then
recast their definition in terms of trees which allows us to give
a new cell decomposition.

\subsection{Spineless Cacti as tree-like ribbon graphs}
We recall from \cite{cact,cyclic} the definition of normalized
spineless cacti and spineless cacti in terms of tree like ribbon
graphs.

\subsubsection{Ribbon graphs}

A ribbon graph is  a graph $\G$ together with a cyclic order
$\prec_v$ of each of the sets $F(v)$. On a ribbon graph, there is
a natural map $N$ from flags to flags given by associating to a
flag $f$ the flag following $\imath(f)$ in the cyclic order
$\prec_{\del(\imath(f))}$. The orbits of the map $N$ are called
the cycles. We say that a vertex lies on a cycle, if there is a
flag in the cycle which is adjacent to this vertex.

If a graph is a ribbon graph, the knowledge of the map $N$ is
equivalent to the knowledge of the cyclic orders $\prec_v$, since
the successor of a flag $f$ is given by $N(\imath(f))$.

Every ribbon graph has a genus which is defined by
$2-2g=\#\mbox{vertices}-\#\mbox{edges}+\#\mbox{cycles}$.\footnote{It
is well known that a ribbon graph can be fattened to a surface
with boundary, such that it is the spine of this surface and the
cycles correspond to the boundary components. The genus is the
genus of the corresponding surface without boundary obtained by
contracting the boundaries to points (or equivalently gluing in
discs).}

A {\em metric} on a graph is a function $\mu:E(\G)\rightarrow
{\mathbb R}_{\geq0}$. A graph with a metric is called a metric
graph.

\begin{df}
A {marked spineless treelike ribbon graph} is a ribbon graph
of genus $0$ together with  a distinguished flag $f_0$,
such that
\begin{itemize}
\item [i)] if $c_0$ is the unique cycle that contains the flag $f_0$, then for each flag
 $f$ either $ f\in c_0$ or $\imath(f)\in c_0$ and
\item[ii)] if $v_0=\delta(f_0)$, then $|F(v_0)|\geq 2$ and for
$v\neq v_0:|F(v)|\geq3$.
\end{itemize}

We will fix the notation $c_0$ for the cycle containing $f_0$ and
call it the distinguished cycle. Likewise we fix the notation
$v_0:=\delta(f_0)$ and call it the root or global zero.
\end{df}

\subsubsection{Spineless cacti}
Cacti without spines are the set of metric marked treelike ribbon graphs.

\begin{df}
We define $\Cact(n)$ to be the set of metric marked treelike
ribbon graphs which have $n+1$ cycles together with a labelling
 of the cycles by $\{0,\dots,n\}$ such that $c_0$ is labelled by $0$.
 We call the elements of $\Cact(n)$ spineless cacti
 with $n$-lobes.
\end{df}

\subsubsection{$\Sn$-action}
Notice that  $\Sn$  acts via permuting the labels $\{1,\dots,n\}$.

\subsubsection{Scaling}
There is an action of ${\mathbb R}_{>0}$  on $\Cact(n)$ which
simply scales the metric. That is if $\mu$ is the metric of the
cactus, the metric scaled by $\lambda\in {\mathbb R}_{>0}$ is
defined by $(\lambda\mu)(e)=\lambda \, \mu(e)$.

\subsubsection{Cactus Terminology} Since the notion of cacti comes with a
history, we set up the usual terminology that is  used in the
literature to describe these objets. Given a spineless cactus with
$n$-lobes, we use the alternate name ``arc'' for ``edge''  and
call $v_0$ the root. Also, we will use the terminology ``special
points'' for the vertices and call the vertices $v$ with
$|F(v)|\geq3$ the intersection points. Sticking with this theme
the arc length of an arc $e$ of a metric spineless treelike ribbon
graph will be simply $\mu(e)$ where $\mu$ is the metric. The lobes
will be the cycles labelled by $\{1,\dots,n\}$ and  the cycle
$c_0$ will be called the perimeter or outside circle. The length
or radius of a lobe or the perimeter is the sum of the lengths of
the underlying edges of the oriented edged belonging to the cycle.

\subsection{Lobes and base points}

Notice that in a spineless cactus, every cycle has a distinguished
flag. For this enumerate the flags of $c_0$ starting with $f_0$.
Now for each cycle $c_i$ other that $c_0$, there is a first flag
$f$ of this linearly ordered set $c_0$ such that $\imath(f)$ is an
element of the cycle $c_i$. The distinguished flag is then defined
to be flag  $N(\imath(f))$. Notice that this flag and $f$ share
the same vertex, the distinguished flag of the cycle. In cactus
terminology these vertices are called the local zeros and
$v_0=\delta(f_0)$ is sometimes referred to as the global zero.

 Also notice that the cycle $c_0$ gives a linear order to all the
edges of a cactus. This linear order is given by fixing the edge $\{f_0,\imath(f_0)\}$ to be the
first edge. In the same fashion, there is a linear order on all edges belonging to the other cycles by
letting the first edge be the one containing the distinguished flag discussed above.

\subsubsection{Dual black and white graph}
Given a marked spineless treelike metric ribbon graph, we can
associate to it a dual graph. This is a graph with two types of
vertices, white and black. The first set of vertices is given by
replacing each cycle except $c_0$ by a white vertex. The second
set of vertices, the set of black vertices, is given by the
vertices of the original graph. The set flags of the dual graph
 is taken to be equal to the set of flags of the original graph.
The maps $\imath$ and $\delta$ are defined in such a way,  that
the edges run only from white to black vertices, where we  two
such vertices are joined if the vertex of the ribbon graph
corresponding to the black vertex lies on the cycle represented by
the white vertex. See the Appendix for the precise combinatorial
description of this construction.

Notice that this is a planted planar bipartite tree.  It is a
planar tree, since the underlying ribbon graph had genus zero, and
the pinning is the one induced by the cyclic order of the flags of
the original graph  (see Appendix). It is planted by defining the
linear order at $v_0:=\del(f_0)$ by letting $f_0$ be the first
flag. Moreover, the set of edges of this tree is in 1--1
correspondence with the set of edges of the ribbon graph and the
two enumerations of the edges by the distinguished cycles agree.
For details, see the Appendix.

\begin{df}
\label{typedef} The topological type of a spineless cactus in
$\cact^1(n)$ is defined  to be the tree $\t \in \wlbptree(n)$
which is its dual b/w graph together with the labelling induced
from the labels of the cactus and the linear order induced on the
edges by the embedding into the plane and the position of the
root.
\end{df}

\begin{figure}
\epsfxsize = \textwidth
\epsfbox{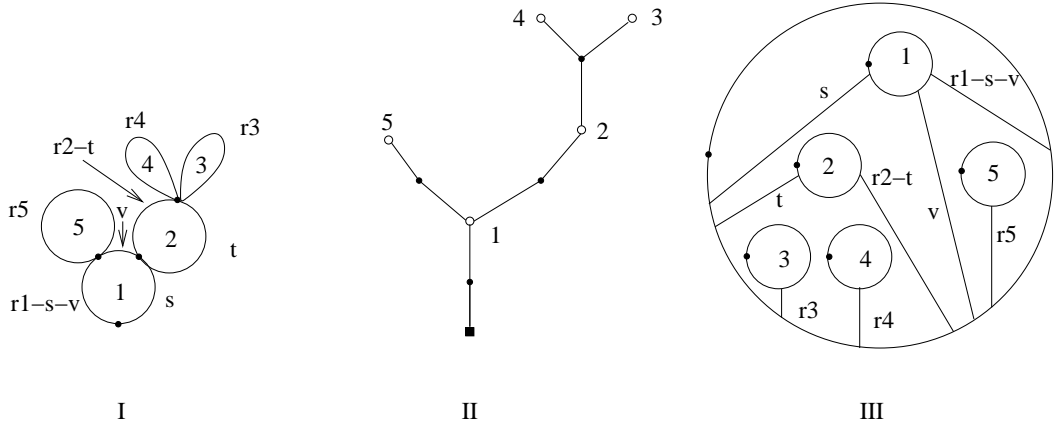} \caption{\label{nospinesdel}I. A cactus without spines with
its arc lengths. II. Its topological type as a labelled planted planar bipartite tree. III. Its depiction as a surface
with boundary and arcs in $\DArc$.}
\end{figure}

For an example of a spineless cactus and its topological type, see Figure \ref{nospinesdel}.

\subsection{The marked treelike spineless ribbon graph obtained from a tree}
\label{arctoedge}

The inverse construction to the dual graph is combinatorially a
little tricky and it can be found in the Appendix.

A good geometric picture is the following. Given a b/w planar
planted bi-partite tree $\t$. Draw the tree in the plane, where we
do not draw a root edge, but rather mark the first flag at the
root vertex. The first part of the construction is to expand each
white vertex to a circle. This can be done, since we have a cyclic
order at each edge. This graph has twice the number of edges as
$\t$, the new edges and the edges of $\t$.  It is a genus zero
ribbon graph and the cycle of $f_0$ contains all edges. The new
edges are in 1-1 correspondence with the edges of $\t$, by
assigning to each edge of $\t$ the edge which is the next edge in
the cyclic order given by $c_0$. To obtain the marked ribbon graph
contract all the edges which were formerly the edges of $\t$. We
choose a marking by fixing the flag which is the successor of the
image of the originally marked flag in the unique cycle of this
graph. The result of this operation is a marked spineless ribbon
graph whose topological type is that of $\t$ and whose edges are
in 1-1 correspondence with those of $\t$. The precise graph
theoretic proof of these statements can be found in the Appendix.

\begin{lem}
\label{typelemma} A spineless cactus $c\in \cact(n)$ is uniquely
determined by its topological type $\t\in \wlbptree(n)$ (as
defined in Definition \ref{typedef}) and the lengths of its arcs.
Moreover the arcs are in 1--1 correspondence to the edges of it
topological type.
\end{lem}

\begin{proof}
By the above results, it follows that each spineless cactus gives
rise to the corresponding data. Vice-versa given $\t\in
\wlbptree(n)$, it gives rise to a marked spineless ribbon graph
 with both constructions being inverses. Since the set of edges of both the ribbon
graph and the tree can be identified, the claim follows.
\end{proof}

\subsubsection{Lobes as circles}
\label{lobeparameters} To simplify the exposition, we now explain
how to think of a cycle as a parameterized circle. Each cycle  is
a sequence of arcs with a given length. So by simply combining
these arcs and viewing the result as a CW complex, we obtain a CW
decomposition of $S^1$ and a metric on this space. In effect the
radius of  this metric $S^1$ is the radius of the lobe or the
outside circle. Moreover there is a distinguished point on this
$S^1$ given by the vertex of the distinguished flag and there is
an orientation induced by the order of the cycle. Thus we can
think of a cycle as a parameterized $S^1$ of radius $r$ where $r$
is the length of the cycle. The technical purely combinatorial
version of this construction is relegated to the Appendix.

A realization of a cactus with $n$-lobes as a CW complex can then
be thought of as a topological space with $n+1$ maps $\phi_i$ from
$S^1_{r_i}$ to this space, so that the $i$-th map is bijective
onto the $i$-th lobe for $i>0$ and $r_i$ is the radius of the
$i$-th lobe. For $i=0$ the map is surjective. In this way we can
think of a cactus in $\cact(n)$ as a collection of  $S^1_{r_i}$'s
which intersect in a tree-like fashion.

Moreover, removing the special points, we are left with a
1-manifold with many components, one for each edge, which is
naturally a subset of the realization of the cactus. We can think
of a lobe as the closure of its components which presents it as a
metric space.

\subsubsection{Gluing for cacti without spines} In this paragraph, we
recall how to define  the operations
\begin{equation}
\label{cactgluing}
 \circ_i: \Cact(n) \times \Cact(m) \rightarrow
\Cact(n+m-1)
\end{equation}
In plain words: given two cacti without spines let
$r_i$ be the length of the $i$-th lobe of the first cactus. First
scale, such that the outside circle of the second cactus has
length $r_i$, then glue in the second cactus by identifying the
outside circle of the second cactus with the i-th circle of the
first cactus.

This gluing is  naturally understood as a gluing of the
CW-complexes. For this realize both graphs, then consider the
$i$-th lobe and the outside circle as a $S^1$s as in
\S\ref{lobeparameters} and then glue the two CW-complexes by
identifying points using the maps of $S^1$  to both spaces as
gluing data. The topologically most satisfying way to define the
gluing is via the Arc operad as it is done in \cite{KLP}. A
detailed recollection of these results would be too much of a
detour, however, to be self-contained, we give the somewhat
technical combinatorial description corresponding to the graph
theoretic definition of spineless cacti in the Appendix.

It is easily checked that the operations (\ref{cactgluing}) turn
$\cact=\{\cact(n)\}$ into an operad \cite{KLP,cact}, albeit for
the moment an operad of sets. We will deal with the topology
shortly.

\subsubsection{Normalized cacti without spines}
Normalized cacti without spines  are the subset of spineless cacti
whose lobes all have length one.

\begin{df} The set of normalized cacti without spines with $n$
lobes is the subset $\cact^1(n)\subset \cact(n)$ of spineless
cacti with $n$-lobes all of whose length is 1.\end{df}

\vskip 3mm

\subsubsection{Gluing for normalized cacti without spines}
\label{glue} We define the operations
\begin{equation}
\circ_i: \Cact^1(n) \times \Cact^1(m) \rightarrow \Cact^1(n+m-1)
\end{equation}
basically by the following procedure: given two normalized cacti without
spines we re-parameterize the $i$-th component circle of the first
cactus to have length $m$ and glue in the second cactus by
identifying the outside circle of the second cactus with the
$i$-th circle of the first cactus. The rigorous combinatorial gluing operation is in the Appendix.

Notice that this gluing differs from the one for spineless cacti,
since now  just a lobe and not a whole cactus is re-scaled.

These maps do not endow the normalized spineless cacti with the
structure of an operad, since they are not associative. But they
induce the slightly weaker structure of a quasi-operad.

\begin{df}\cite{cact}
\label{compositordef} A  quasi-operad is an operad where the
associativity need not hold. An operad in the category of
topological spaces is called homotopy associative if the
compositions are associative up to homotopy.
\end{df}

\subsection{Scaled cacti and other variations}

\subsubsection{Projective Cacti}
\label{cactscale} Spineless cacti come with a universal scaling
operation of ${\mathbb R}_{>0}$ which simultaneously scales all
radii by the same factor $\lambda \in {\mathbb R}_{>0}$. This
action is a free action and the gluing descends to the quotient by
this action. We call the resulting operad projective spineless
cacti and denote its spaces by $\Pcact(n) :=
\cact(n)/\mathbb{R}_{>0}$.

It is clear that $\cact(n)=\Pcact(n)\times \mathbb{R}_{>0}$.

\begin{rmk}
Notice that the topological type as defined in Definition
\ref{typedef} is invariant under the global re-scaling action of
$\mathbb{R}_{>0}$ so that a projective spineless cactus also has a
well defined topological type.
\end{rmk}

\subsubsection{Left, right and symmetric cacti operads}
\label{leftrightsymmetric} For the gluing maps (\ref{cactgluing})
one has three basic possibilities to scale in order to make the
size of the outer loop of the cactus that is to be inserted match
the size of the lobe into which the insertion should be made.
\begin{enumerate}
\item Scale down the cactus which is to be inserted. This is the
original version - we call it the right scaling version. \item
Scale up the cactus into  which will be inserted. We call it the
left scaling version. \item Scale both cacti before gluing. Let
$r_i$ be the length of the lobe $i$ of $c$ and let $R=\sum_j r_j'$
be the length of the outside circle of $c'$. Now to define
$c\circ_i c'$ first scale $c$ by $R$ and $c'$ by $r_i$. Then
identify the outside loop of $c'$ with the lobe $i$ of $c$ which
now both have length $Rr_i$. We call this the symmetric scaling
version.
\end{enumerate}

All of these versions turn out to be homotopy equivalent. In
particular passing to the quotient operad $\Pcact$ they all
descend to the same operations.

The advantages of the different versions are as follows: version
(1) is the original one and inspired by the re-scaling of loops,
i.e.\ the size of the outer loop of the first cactus is constant.
Version (2) has the advantage that cacti whose lobes have integer
sizes are a suboperad. And version (3) is the one which needs to
be used to embed  into the operad $\darc$ of \cite{KLP}. In this
version there is an embedding of the operad of spineless cacti
operad into the cyclic operad $\darc$. Projective spineless cacti
embed into the cyclic operad $\Arc$. The equivalent statements for
the larger operad of cacti which contains the operad of spineless
cacti under consideration as a suboperad also hold true.

\subsection{The topology} In this section, we give a short
account of how to put a topology on the set of spineless cacti. As
discussed above, the quickest way would be to give the topology to
the spaces $\cact^1$ as subspaces of the operad $\DArc$ as defined
in \cite{KLP}, see also \cite{cact}. We will also give an
equivalent way to define the topology in the next section by
identifying normalized spineless cacti with a CW complex. This
description is intrinsic and the most adapted to cacti. This is
the reason why we decided to define the topology for spineless
cacti and cacti using the CW complex approach of this article in
\cite{cact} rather than using one of the other equivalent
descriptions. Endowed with a topology spineless cacti become a
topological operad and normalized spineless cacti a homotopy
associative topological quasi-operad.

\subsubsection{Cacti glued from open cells}
For definiteness, we give one construction of the topology which
is tantamount to defining the topology of $\Pcact$ as an open
subset of $\Arc$ without referring the reader to \cite{KLP}.

\begin{nota}
Let $\Delta^n$ denote the standard $n$-simplex and $\dot \Delta^n$
its interior.

For $\t\in \wlbptree$, set $\SE(\t)=\D^{|E_w(\t)|}$ and denote its interior by
$\dot{\SE(\t)}$.
\end{nota}

\begin{lem}
The projective spineless cacti of a fixed topological type $\t$ are in
bijection with points of the interior $\dot{\SE(\t)}$ of the
simplex $\SE(\t)$. Moreover
$$
\Pcact(n)=\amalg_{\t\in \wlbptree(n)}\dot{\SE(\t)}
$$
\end{lem}

\begin{proof}
First notice that in each class of a projective spineless cactus
there is a unique representative
 whose arc lengths sum up to
one.   Using barycentric co-ordinates on the simplex $\SE(\t)$ we
can thus  identify the projective spineless cacti of the given
topological type $\t$ with points in the open simplex
$\dot{\SE(\t)}$. Hence the first bijection follows from Lemma
\ref{typelemma}. It is clear that every topological type occurs,
so the second statement follows.
\end{proof}

\begin{df}
We define the degeneration of a spineless cactus $c$ with respect
to an arc $a$ which is not an entire lobe to be the spineless
cactus obtained by contracting the arc $a$. If the root was on the
boundary of $a$, the image of $a$ is the new root. The marked lobe
is the image of the lobe to which the arc immediately following
$a$ around the outside circle belonged.
\end{df}

\begin{rmk}
Notice that if $\t$ is the topological type of a spineless cactus
$c$ and $e=\{w,b\}$ is the edge corresponding to the arc $a$ in
the terminology of Remark \ref{arctoedge} then the topological
type of the degeneration of the spineless cactus $c$ with respect
to the arc $a$ is $\del_e(\t)$.
\end{rmk}

\begin{df}
\label{topdef} We give $\Pcact(n)$ the topology induced by
identifying the $e$-th open face of $\SE(\t)$ with
$\dot{\SE(\del_e(\t))}$ for any $e\in E_{ angle}(\t)$.

We define the topology for $\cact(n)=\Pcact(n)\times
\mathbb{R}_{>0}$ to be the product topology and endow
$\cact^1(n)\subset \cact(n)$ with the subset topology.
\end{df}

This identifies the $e$-th open face of $\SE(\t)$ with the
degenerations of spineless cacti in $\SE(\t)$ with respect to the
arc $a$ that corresponds to $e$.

\subsubsection{Lemma}
\label{toplemma}
 As topological spaces $\Cact(n)=\Cact^1(n)\times
\mathbb{R}^n_{>0}$ and $\Pcact(n)=\Cact(n)/\mathbb{R}_{>0}$.
\begin{proof}
The first statement follows by identifying the
 factors $\mathbb{R}^n_{>0}$ with the sizes of the lobes.
 The second statement follows from the observation that
 the $\mathbb{R}_{>0}$ action is free and continuous.
\end{proof}

\begin{prop} The gluing maps (\ref{cactgluing}) endow the spaces
$\Cact(n)$ with the structure of a topological operad.
\end{prop}
\begin{proof}
Straightforward.
\end{proof}

\begin{prop}\cite{cact}
The glueings of \ref{glue} together with the permutation action of
$\Sn$ on $\Cact^1(n)$ turn $\Cact^1(n)$ into a topological
quasi-operad which is homotopy associative.
\end{prop}
\begin{proof}
Straightforward.
\end{proof}

\subsection{Other approaches to the topology}
In this subsection we give a brief non-technical summary of
alternative approaches to give spineless cacti a topology by
identifying them with other well known objects.  Although the
approaches seem vastly different they all lead to the same
topology. For a detailed but still concise summary see
\cite{cact}.

\subsubsection{Topology of Cacti as special types of ribbon graphs}
A topology on this set  of cacti could alternatively be given as follows:
given by the metric and the following
convention. If the length of an arc goes to zero, the respective
edge is contracted. If one of the two vertices of the edge
corresponded to a root, the new root marking will be the vertex
corresponding to the contracted edge. The cycle of the root
marking is defined as the image of the cycle to which the edge
immediately preceding the contracted edge in the order of the
distinguished cycle belonged. Vice-versa the realization of a graph as
above defines a unique normalized spineless cactus and hence the
set of normalized spineless cacti is topologized.

\subsubsection{Cacti-Ribbon graphs as surfaces with weighted arcs}
To get an element of $\DArc$ take the ribbon graph and
consider a surface of which it is the spine. This surface will be
of genus zero and can be realized as a disc with holes. Run arcs
from the inside boundaries (viz.\ the ``holes'') to the boundary,
one for each edge, such that it crosses this edge transversally
and has no other crossings with edges or arcs. In other words the
ribbon graph is the dual graph on the surface to the graph given
by the arcs and the boundaries.  Mark a point on each boundary,
which is not the endpoint of an arc, so that the linear order on
the cycles of the ribbon graph agrees with the linear order on the
arcs induced by
 going around the boundary in the orientation of the disc starting at
the marked point. Each arc carries a weight given by the value of
the metric on the respective edge to it. This gives the unique
element of $\DArc$ representing the spineless cactus. In the limit
where a weight goes to zero, the respective arc is simply erased
\cite{KLP}.

For an example, see Figure \ref{nospinesdel}.

\subsection{The relation between spineless cacti and normalized
spineless cacti}

In this section we briefly review the relationship between
spineless cacti and their normalized version. For the full
details, we refer to \cite{cact}.

\subsubsection{The scaling operad}
We define the scaling operad  ${\mathcal R}_{> 0}$ to be given by
the spaces ${\mathcal R}_{> 0}(n):=  \mathbb {R}_{> 0}^{n}$ with
the permutation action by $\Sn$ and the following products

$$
(r_1,\dots,r_n) \circ_i (r'_1,\dots,r'_m) = (r_1, \dots r_{i-1},
\frac{r_i}{R} r'_1, \dots, \frac{r_i}{R} r'_m, r_{i+1},\dots
r_{n})
$$
where $R=\sum_{k=1}^m r'_k$.

On one hand the scaling operad keeps track of the sizes of the
lobes of a cactus. On the other hand the difference between the
compositions in normalized and non-normalized spineless cacti is
given by an action of the following type:
\begin{equation}
\rho_i: Cact^1(n) \times {\mathcal R}_{> 0}(m) \times Cact^1(m)
\rightarrow Cact^1(n)
\end{equation}
The effect of the action is to move the intersection points of the
lobes incident to the lobe $i$ of the first normalized spineless
cactus around that lobe according to the outside circle of
spineless cactus to be inserted. The action is defined so that
after the displacement of the lobes the composition of the
perturbed normalized spineless cactus with the other normalized
spineless cactus as normalized spineless cacti coincides with the
composition obtained by gluing the two normalized spineless cacti
simply as spineless cacti and then scaling back each lobe of the
resulting cactus to length one. For the explicit formulas, we
refer to \cite{cact}. Using this action, we can perturb the
multiplications of normalized spineless cacti to fit with those of
spineless cacti.
\begin{multline}
\circ_i^{{\mathcal R}_{>0}}: \Cact^1(n) \times
{\mathcal R}_{> 0}(m) \times  \Cact^1(m)\\
\stackrel{id\times id \times \Delta}{\longrightarrow} \Cact^1(n)
\times {\mathcal R}_{>
0}(m) \times  \Cact^1(m)\times  \Cact^1(m)\\
 \stackrel{\rho_i \times id}{\longrightarrow}
\Cact^1(n) \times \Cact^1(m) \stackrel{\circ_i}{\longrightarrow}
\Cact^1(n+m-1)
\end{multline}
to get a continuous map
$$
(c,\vec{r'},c') \mapsto c \circ_i^{\vec{r'}}c'
$$

Due to the nature of the maps above it is possible to continuously
``un-deform'' the deformed product while staying in the category
of quasi-operads.

\begin{thm}\cite{cact}
\label{semicact} The operad of spineless cacti is isomorphic to
the operad given by the semi-direct product of its normalized
version with the scaling operad. The latter is homotopy equivalent
(through quasi-operads)  to the direct product as a quasi-operad.
The direct product is in turn equivalent as a quasi-operad to
$\Cact^1$.

\begin{eqnarray}
\Cact &\cong&   {\mathcal R}_{>0}\ltimes \Cact^1
\sim  \Cact^1 \times {\mathcal R}_{>0}\simeq \Cact^1
\end{eqnarray}
where the semi-direct product compositions are given by
\begin{equation}
(\vec{r},c)\circ_i(\vec{r}',c') = ( \vec{r}\circ_i \vec{r'}, c
\circ_i^{\vec{r'}} c')
\end{equation}
\end{thm}

From this description one obtains several useful corollaries
\cite{cact}. The ones relevant to the present discussion are listed below.

\begin{cor}
The quasi-operad of normalized spineless cacti is homotopy
associative and thus its homology quasi-operad is an operad.

Moreover as quasi-operads normalized spineless cacti and spineless
cacti are homotopy equivalent via a homotopy of quasi-operads.
\end{cor}

And finally:

\begin{cor}
\label{quasiiso} Normalized spineless cacti are operadically
quasi-isomorphic spineless cacti. I.e.\ their homology operads are
isomorphic.
\end{cor}

\subsection{Spineless cacti and the little discs operad}

The most important result of \cite{cact} which we will use is:

\begin{thm}\cite{cact}
\label{littlediscscacti} The operad $\cact$ is equivalent
to the little discs operad.
\end{thm}

\section{A cell decomposition for normalized spineless cacti}

\subsection{The cell complex}

\begin{rmk}

\label{lengthrem} For a normalized spineless cactus the lengths of
the arcs have to sum up to the radius of the lobe and the number
of arcs on a given lobe  represented by a vertex $v$ is $|v|+1$.
Hence the lengths of the arcs lying on the lobe represented by a
vertex $v$ are in 1-1 correspondence with points of the simplex
$|\Delta^{|v|}|$. The coordinates of $|\Delta^{|v|}|$ naturally
correspond to the arcs of the lobe $v$ on one hand and on the
other hand to the incident edges to $v$ in the dual b/w graph.
\end{rmk}

\begin{df}
We define $\wlbptree(n)^k$ to be the elements of $\wlbptree(n)$
with $|E_w|=k$.
\end{df}

\begin{df} For $\t \in \wlbptree$ we define
\begin{equation}
\D(\t):=\times_{v \in V_w(\tau)}\D^{|v|}
\end{equation}
We define $C(\t)=|\D(\t)|$.  Notice that $\dim(C(\t))=|E_w(\t)|$.

Given $\D(\t)$ and a vertex $x$ of any of the constituting
simplices of $\D(\t)$ we define the $x$-th face of $C(\t)$ to be
the subset of $|\D(\t)|$ whose points have the $x$-th coordinate
equal to zero.
\end{df}

\begin{df}
We let $\CWCact(n)$ be the CW complex whose k-cells are indexed by
$\t \in \wlbptree(n)^k$ with the cell $C(\t)=|\D(\t)|$ and the
attaching maps $e_{\t}$ defined as follows. We identify the $x$-th
face of $C(\t)$ with $C(\t')$ where $\t'=\del_e(\t)$ is the
topological type of the  spineless cactus $c'$ which is the
degeneration of a spineless cactus $c$ of topological type $\t$
with respect to the arc $a$ that simultaneously represents the
vertex $x$ of $\D(\t)$ and the edge $e$ of $\t$ (see Remark
\ref{arctoedge}).
\end{df}

We denote by $\dot e_{\t}$ the restriction of $e_{\t}$ to the
interior of $\D(\t)$. Notice that $\dot e_{\tau}$ is a bijection.

\begin{thm} \label{cact1} The space $\Cact^1(n)$ is homeomorphic to
 the CW complex $\CWCact(n)$.
\end{thm}

\begin{proof}
The map from $\Cact^1(n)$ to $\CWCact(n)$ is given by Lemma
\ref{typelemma} and Remark \ref{lengthrem}. Vice versa given an
element on the right hand side, the unique open cell it belongs to
determines the topological type and it is obvious that any tree in
$\wlbptree$ can be realized.  Then the barycentric coordinates
assign weights to the arcs via the correspondence of the vertices
of the factors of $\D(\t)$ with arcs of the cactus.

For the homeomorphism, we notice that the bijection restricted to
the insides of the top dimensional cells is obviously a
homeomorphism. This is seen by slightly perturbing the non-zero
lengths of the arcs. The limit where  one of the lengths of the
arcs goes to zero is given by passing to the corresponding face of
the corresponding simplex factor of $C(\t)$. The resulting tree,
which is the topological type of the limit cactus, will be the
tree which was used in the definition of the attaching map. Thus
the contraction of the arc corresponds to the projection map
sending the respective coordinate to zero. This agrees with the
limit in $\SE(\t)$. Therefore in the case of a degeneration the
limits also agree and the result follows.
\end{proof}

\subsubsection{Chains for $\cact$}
\label{cactcells}
 Since the factors of $\mathbb{R}_{>0}$ are
contractible,  it is clear that $\CCcact$ is a chain model for
$\cact$. Furthermore, by  Theorem \ref{celldecomp} below, we will
see that $\CCcact$ is even an operadic chain model for $\cact$ and
hence by Theorem \ref{littlediscscacti} for the little discs
operad.

Instead of using purely the cells of $\cact^1$ as a chain model
for the little discs operad one can choose any operadic chain
model $Chain(\mathcal {R}_{>0})$ for the scaling operad and then use
the mixed chains for $\cact$ i.e.\ $\CCcact\otimes
Chain(\mathcal{R}_{>0})$. It follows that the inclusion of the
cellular chains of $\cact^1$ into the mixed chains is an inclusion
of operads up to homotopy.

Given an operation of the $\CCcact$ we can let the mixed chains of
$\cact$ act by letting the action of the mixed chains of bi-degree
$(n,0)$ be that of the component of $\cact^1$ and setting the
action of all the other chains to zero.

In any case, $\CCcact$ is chain equivalent to any form of chain
complex of mixed chains $\CCcact\otimes Chain(\mathcal{R}_{>0})$.

\subsubsection{Pseudo-Cells for $\Pcact$}
From Theorem \ref{cact1} and Lemma \ref{toplemma}, we obtain a
decomposition for $\Pcact$ as
$$\Pcact(n)=\amalg_{\t \in \wlbptree(n)} \tilde C(\t), \quad \tilde C(\t)
= \D(\t) \times \dot\Delta^n.$$

This is not a cell decomposition since this $\tilde C(\t)$ are not
really cells, but this decomposition will be useful in the following.

\begin{rmk} One could formally close the $\tilde C(\t)$ to cells
and glue them by forgetting the lobes whose lengths goes to zero
as described. This type of forgetting map is a quasi-fibration
as shown in \cite{cact}. The gluing would yield a CW space
which is the union over $n$ of all $\Pcact(n)$.
We will not pursue this construction further here.
\end{rmk}

\subsection{Orientations of Chains}
\label{genfix}
To fix the generators and thereby the signs for the chain operad
we have several choices, each of which is natural and has appeared
in the literature.

To fix a generator  $g(\t)$ of $\CCcact$  corresponding to the
cell indexed by $\t\in \wlbptree(n)$ we need to specify an
 orientation for it, i.e.\ an explicit parametrization or
 equivalently an order of the white edges of the tree it is represented by.

The first orientation which we call $\nat$ is the orientation
given by the natural orientation for a planar planted tree. I.e.\
fixing the order of the white edges to be $\prec^{\t}$. N.B.\ this
actually coincides with the natural orientation of the cells of
$\Arc$, see \cite{KLP}.

We will also consider the orientation $\Op$ which is the
enumeration of the white edges which is obtained by starting with
the incoming edges of the white vertex labelled one, in the order
$\prec^{\t}_{v_1}$, then continuing with the incoming white edges
of the vertex two, etc.\ until the last label is reached.

Finally, for top-dimensional cells, we will consider the
orientation of the edges induced by the labels, which we call
$\lab$. It is obtained from $\nat$ as follows: for $\t \in
\wptree$ let $\s\in \Sn$ be the permutation which permutes the
vertices $v_1, \dots, v_n$ to their natural order induced by the
order $\prec^{\t}$. Then let the enumeration of $E_w$ be
$\s(\nat)$, where the action of $\s$ on $E_w$ is given by the
correspondence $\verttoedge$ and the correspondence between black
and white edges via $(v,N(v))\mapsto (N(v), N^2(v))$ for top
dimensional cells.

To compare with the literature it is also useful to introduce the
orientations $\bnat$, $\blab$, and $\bop$ which are the reversed
orientation of $\nat$, $\lab$ and $\Op$, i.e.\ reading the
respective orders from right to left.

\begin{lem}
\label{keylemma}
In the orientation $\nat$ the induced quasi-operad structure on the
cellular chains of $\CCcact$ is given by
$$
C(\t)\circ_i C(\t')=C(\t\circ_i \t')
$$
where we understand the right hand side to be given by the natural
extension to  $\mathbb{Z}$-modules, i.e.\ $C(\sum_i z_i\t):=\sum_i
z_i C(\t)$.
\end{lem}

\begin{proof}
Let $\t\in \wlbptree(n)$ and $\t'\in \wlbptree(m)$. If we glue a
cactus from $C(\t')$ into $C(\t')$ at the lobe $i$, its
topological type will be one of the trees in the sum
$\t\circ_i\t'$. Conversely, fix a tree $\tilde \t$ which appears
in the sum $\t \circ_i \t'$. Now any element $\tilde c$ of
$C(\tilde \t)$  can be uniquely decomposed as $c \circ_i c'$ with
$c \in C(\t)$ and $c' \in C(\t')$ as follows. Since the
topological type $\tilde \t$ is one of the summands of $\t \circ_i
\t'$ if follows that  the lobes $i$ through $i+n$  are connected.
Let $c'$ be the normalized spineless sub-cactus which consists of
the lobes $i$ through $i+m-1$ and has  the local zero of the lobe
$i$ as the root marking. The cactus $c$ will be the cactus
constituted by the lobes $1$ through $i-1$ and $i+n$ through
$n+m-1$ together with a lobe marked $i$ which is the outside
circle of $c'$ re-parameterized to the length one. Its root will
be the root of $\tilde c$ if it does not lie on $c'$. And if it
does, it will necessarily be the local zero of the $i$-th lobe of
$\tilde c$ and the root of $c$ will be this point thought of as
the zero of the outside circle of $c'$ which is the $i$-th lobe of
$c$.

 To give the co-ordinates of this construction let $\t'$ be
the connected subtree of $\tilde\t$ whose white vertices are the
white vertices labelled $i$ through $i+m-1$ and whose black
vertices are the $N(v_j),j\in\{i,\dots,i+m-1\}$ together with the
induced structure of planar planted tree. Let $\t:=\tilde \t/\t'$
be the planar planted tree obtained from $\tilde \t$ by
contracting $\t'$. Let the coordinates of $\tilde c$ in $C(\tilde
\t)$ be $v=(\mathbf {v}_{1},\dots \mathbf{v}_{n+m-1})$, were we
use the short hand notation
$\mathbf{v}_{i}=({v}_{v_i,1},\dots,{v}_{v_i,|v_i|+1})\in
|\D^{|v_i|}|$. Then $c'$ is the cactus with coordinates
$(\mathbf{v}_i,\dots, \mathbf{v}_{n+m-1})$ in $C(\t')$. To define
$c$ let $(d_1,\dots,d_k)$ be the non-zero distances along the
outside circle between the lobes of $\tilde \t$ meeting the
subtree $\t'$ which are not part of $\t'$. Now $c$ is the cactus
with coordinates $v=(\mathbf {v}_{1},
\mathbf{v}_{{i-1}},\mathbf{\bar v}, \mathbf{v}_{{i+m-1}},\dots
\mathbf{v}_{n+m-1})$ in $C(\t)$ where
$\mathbf{\bar{v}}=\frac{1}{m}(d_1,\dots,d_k)$. It is clear that
$(c,c')$  will be the only pre-image of $\tilde c$. Hence the sets
are in bijection and therefore the sum contains the correct
summands with the correct multiplicity up to a sign.

If one wishes to construct a representing configuration one just
has to resolve the intersection points of $c'$. This can be done
by for instance choosing a thickening of the cactus to a surface
taking the boundary corresponding to the outside circle and
attaching the lobes  at the distances around this circle as
dictated by the cactus $\tilde c$.

To verify the signs, we have to check that the respective
orientations agree. For this we notice that the coordinates in the
orientation $\nat$ correspond to the  white edges enumerated in
the natural order of the given tree  and that the permutation
induced on the co-ordinates is exactly the sign incorporated into
the definition of $\t\circ_i \t'$.  Thus the signs agree as
stated.
  \end{proof}

\begin{lem}
On the level of sets let $\tilde C(\t\circ_i^+
\t')=\amalg_{\tilde\t} \tilde C(\tilde \t)$ where $ \tilde \t$
 is a summand of $\t\circ_i^+\t'$ then the map
 $\circ_i:\tilde C(\t)\times \tilde C(\t') \rightarrow \tilde C(\t\circ_i^+ \t')$ is a bijection.
\end{lem}

\begin{proof}
As above, for the equation  $\tilde c= c\circ_i c'$,  identify
$c'$ and $c$ with the respective spineless sub-cactus and quotient
cactus, yielding a unique pre-image.
\end{proof}

This gives a topological proof for the associativity of the
$\circ_i^+$.

\begin{cor}
\label{topasscor} The operations $\circ_i^+$ give $\wlbptree$ the
structure of an operad.
\end{cor}

\begin{prop}
\label{treescelliso} For the choice of orientation $\nat$ and the
induced operad structure $\circ_i$ the map $\t \mapsto g(\t)$
where $g(\t)$ is the generator corresponding to $C(\t)$ fixed in
\S \ref{genfix} is a map of differential graded operads which
identifies $\wlbptree(n)^k$ with $\CCkcactn$, where $CC_k$ are the
dimension-$k$ cellular chains.

The same holds true for the orientation $\Op$ with the appropriate
changes to the signs of the operad $\wlbptree$ discussed in
\S\ref{signchoices}. Finally the analogous statement holds true
when passing to operads indexed by sets for both $\cact$ and
$\wlbptree$.
\end{prop}

\begin{proof}
This statement for the orientation $\nat$ follows from the Lemma \ref{keylemma} above.
For the other orientations the signs will agree, since they are forced onto
the combinatorial operad by definition.
\end{proof}

\begin{thm}
\label{celldecomp} The glueings induced from the glueings of
spineless normalized cacti make the spaces $\CCcactn$ into a chain
operad. Thus $\CCcact$ is an operadic model for the chains of the
little discs operad.
\end{thm}

\begin{proof}
The first statement follows from Proposition \ref{treeopprop}
together with the Proposition \ref{treescelliso}. The second part
of the statement follows from Theorem \ref{littlediscscacti} and
 Corollary \ref{quasiiso}.
\end{proof}

\subsection{Operadic action of $\wlbptree$}
\label{operadicaction} Given an assignment of operations on a
complex $(\O,\delta)$ to the trees $l_n$ and $\t_n^b$ of Figure
\ref{nleaf}, a natural way to let a tree $\t\in \wlbptree$ act on
the complex $(\O,\delta)$ is given as follows. Formally read off
the operation of $\t$ by decorating the white vertices with
elements of $\O$ and interpreting the tree as a flow chart,
assigning to each white vertex $w$ the operation of $l_{|w|}$ and
to each black vertex $v$ the operation of $\t_{|v|}^b$. If the
operations of $l_n$ and $\t_n^b$ are also compatible with the
differentials, then one can obtain a dg-action of $\wlbptree$ by
regarding a mixed complex $CC_*(\cact^1) \otimes \underline
{Hom}(\O,\O)$ in which the order of the tensor factors of
homogenous elements of the tensor product is dictated by the trees
$\t\in\wlbptree$. For this one identifies the white edges of a
tree with the co-ordinates from $\cact^1$ and the vertices with
elements from $\underline {Hom}(\O,\O)$ and builds a differential
from the tree differential  and the internal differential of $\O$.

In the graded case, a sign for these operations and the
differentials has to be included. A way for implementing this idea
compatibly is as discussed in the following.

\subsubsection{Tensor orders}
 For an action of the operad $\wlbptree$ on a
collection of objects $O(n_i)$  of a monoidal category
$\mathcal{C}$,
 we will consider maps
\begin{eqnarray}
\rho: \wlbptree(k)
&\rightarrow& Hom(O(n_1) \otimes \dots \otimes O(n_k),O(m))
\nn\\
\label{opmaps} \t \mapsto [f_1\otimes  \dots\otimes  f_k &\mapsto&
\t(f_1 \otimes \dots \otimes  f_k)]
\end{eqnarray}
Actually, $\t(f_1 \otimes \dots \otimes f_k)$ will be zero unless
  $m=\sum_{v\in E_w(\t)} (n_{\lb^{-1}(v)}-|v|)$.

In the cases of interest for the present considerations, one can
furthermore associate an object of $\mathcal{C}$ to each $\t\in
\wlbptree$ which is usually of the form $\bigotimes_{e\in E_w(\t)}
L$ with $L$ of dimension or degree plus or minus one. We call this
the decomposable case.

Alternatively one can sometimes associate an object of the form
$\bigotimes_{v\in E_w(\t)} D_{|v|}$ to $\t$ where the dimension or
degree of $D_{|v|}$ is $|v|$, e.g.\ the simplex $\Delta^{|v|}$. We
call this the operadic case.

If this is not possible, one can still consider a mixed complex
made up of tensors of objects $C(\t)$ and the $O(n_i)$.

In all these descriptions if the operad has a differential there
is a natural differential obtained by using a combination of the
tree differential and the operadic differential on the various
factors. We now make these constructions precise.

\subsubsection{Tensor order for the decomposable case}
Fix $\mathcal{C}$ to be $\Chain$ or $\Vect$. In this case we wish
to consider the maps (\ref{opmaps}) formally as maps

\begin{equation}
\label{tensordarstellung}
 \bigotimes_{e\in E_w(\t)}L\otimes\bigotimes_{v\in V_w(\t)}
 O_{n_{\lb^{-1}(v)}}\rightarrow O_m
\end{equation}

In the graded case, we have to fix the order of
(\ref{tensordarstellung}). We do this by using $\prec^{\t}$ to
give the tensor product the natural operadic order. This amounts
to formally inserting $O(n_i)$ into the vertex $v_i$. For this let
$N:= |V_w(\t)\amalg E_w(\t)|$ and again let $\num: V(\t)\amalg
E(\t) \rightarrow \{1,\dots,N\}$ be the bijection which is induced
by $\prec^{\t}$. We fix $L$ to be a ``shifted line" i.e.\ a free
$\mathbb{Z}$-module or $k$ vector space generated by an element of
degree plus one (or minus one). Now set
\begin{equation}
W_i :=
\begin{cases} O(n_j)&\text{ if } \num^{-1}(i) = v_j\\
L&\text{ if $\num^{-1}(i)$ is a white edge}
\end{cases}
\end{equation}

We then define the order on the tensor product on the l.h.s.\ of the
expression (\ref{tensordarstellung}) to be given by
$$
W:= W_1 \otimes \dots \otimes  W_{N}
$$

Another way to add the necessary signs but circumnavigating the
use of maps of the type (\ref{tensordarstellung}) would be to
define the operation $\rho$ to include the sign which is the sign
of the shuffle of the two ordered sets on the l.h.s.\ of
(\ref{tensordarstellung}) ordered by the orders on the edges and
vertices into the order  given by $W$. This has the drawback of
using a new convention for each $\rho$.

\subsubsection{The operadic case}
In case we can associate objects $D_n$ in $\mathcal{C}$ to each
n-ary white vertex, we would consider the maps
(\ref{opmaps}) formally as maps

\begin{equation}
\label{alttensordarstellung}
 \bigotimes_{v\in V_w(\t)} ( O_{n_{\lb^{-1}(v)}}\otimes D_{|v|})\rightarrow O_m
\end{equation}
where the  order for the tensor product on the l.h.s.\ is defined
to be the one given by $\prec^{\t}$ on the white vertices.

\subsubsection{The action of the symmetric group}
The action of the symmetric group on the maps $\rho$ is induced by
permuting the labels and permuting the elements $O_{n_i}$. This
induces  signs by pulling back the permutation onto the factors of
$W$. These signs should be included in the definition of
$\Sn$-equivariance of the operadic action.

\begin{rmk}
This treatment of the signs is essential if one is dealing with
operads vs. non-$\Sigma$ operads and wishes to obtain equivariance with respect to the
symmetric group actions. In general the symmetric group action on
the endomorphism operads will not produce the right signs needed
in the description of the iterations of the universal
concatenation $\circ$ of \S \ref{metaops}. In particular this is
the case for Gerstenhaber's product on the Hochschild cochains.
The above modification however leads to an agreement of signs for
the action of the symmetric group for the subcomplex of the
Hochschild complex generated by products and the brace operations,
see \S \ref{delsection}. Another approach is given by viewing the
operations not as endomorphisms of the Hochschild cochains but
rather as maps of the Hochschild cochains to Hochschild cochains
twisted by tensoring with copies of the dual line $L^*$.

If one is not concerned with the action of the symmetric
group, then one can forgo this step.
\end{rmk}

\subsection{Examples}
An example of the type of action described above as the operadic
case is given in \cite{KLP} by the action of $Chain(\Arc)$ on
itself. For the homotopy Gerstenhaber structure we should consider
the chains $CC_*(\cacti^1)$ and  any choice of chain model for
$\Arc$ or any of the suboperads which are stable under the action
of the linear trees suboperad which is the image of spineless
cacti.

The operadic ordering was also used in \cite{CS1} to define the
action of string topology on the chains of the free loop space of
a compact manifold.

We get agreement with the usual signs and conventions of
Gerstenhaber's original results \cite{G} upgraded to operads (see
\S \ref{metaops}) with those of the operations of $\Arc$ and those
of string topology \cite{CS1}, if we denote the action of $\t_1$
as $*^{op}$ and $\t^b_2$ as $\cdot$; see \cite{KLP} for the
operations and Figure \ref{nleaf} for the definitions of the
trees.

\section{Spineless cacti as a natural solution to Deligne's conjecture}
\label{delsection}

\subsection{The Hochschild complex, its Gerstenhaber structure and
Deligne's conjecture}

Let $A$ be an associative algebra over a field $k$. We define
$CH^*(A,A) := \bigoplus_{q\geq 0}CH^q(A,A)$ with $CH^q(A,A)=
\mathrm{Hom}(A^{\otimes q},A)$.

There are two natural operations
\begin{eqnarray*}
\circ_i:CH^p(A,A)\otimes CH^q(A,A)&\rightarrow& CH^{p+q-1}(A,A)\\
 \cup: CH^n(A,A)\otimes CH^m(A,A)&\rightarrow&CH^{m+n}(A,A)
\end{eqnarray*}
where the first morphism is for $f \in CH^p(A,A)$ and $g\in
CH^q(A,A)$

\begin{equation}
f\circ_ig(x_1,\dots,x_{p+q-1})=f(x_1,\dots,x_{i-1},g(x_i,\dots, x_{i+q-1}),
x_{i+q}, \dots, x_{p+q-1})
\end{equation}
and the
second is given by the multiplication
\begin{equation}
f(a_1\dots,a_m)\cup g(b_1,\dots,b_n) =
f(a_1\dots,a_m)g(b_1,\dots,b_n)
\end{equation}

\subsubsection{The differential on $CH^*$}
The Hochschild complex has a differential which is derived from
the algebra structure.

Given $f\in CH^n(A,A)$ then
\begin{multline}
\del(f)(a_1, \dots , a_{n+1}):=
a_1f(a_2, \dots,a_{n+1}) - f(a_1a_2, \dots, a_{n+1})+ \\
 \cdots + (-1)^{n+1}   f(a_1, \dots, a_n a_{n+1}) + (-1)^{n+2}
f(a_1,\dots, a_n)a_{n+1}
\end{multline}

\begin{df}
The Hochschild complex is the complex $(CH^*,\del)$, its cohomology
is called the Hochschild cohomology and denoted by $HH^*(A,A)$.
\end{df}

\subsubsection{The Gerstenhaber structure}
Gerstenhaber \cite{G} introduced the $\circ$ operations:
for $f\in CH^p(A,A)$ and  $g\in CH^q(A,A)$
\begin{equation}
f\circ g:= \sum_{i=1}^p (-1)^{(i-1)(q+1)} f\circ_i g
\end{equation}
and defined the bracket
\begin{equation}
\{f, g\} := f\circ g - (-1)^{(p-1)(q-1)}g\circ f
\end{equation}
and showed that this indeed induces what is now called a
Gerstenhaber bracket, i.e.\ an odd Poisson bracket for $\cup$, on
$HH^*(A,A)$. Here odd Poisson bracket means odd Lie bracket and
the derivation property of the bracket with shifted (odd) signs.

\subsection{Contents of Deligne's Conjecture}
Since $HH^*(A,A)$ has the structure of a Gerstenhaber algebra one
knows from general theory \cite{cohen,cohen2} that thereby
$HH^*(A,A)$ is an algebra over the homology operad of the little
discs operad $D_2$. Now the Gerstenhaber structure on $HH^*(A,A)$
actually stems from the cochain level and the operadic structure
of $H_*(D_2)$ even originates from the level of topological
operads. The question of Deligne was:
\begin{quest}\cite{letter} { Can one lift the
action of the homology of the little discs operad to the chain
respectively cochain level? Or in other words: is there a chain
model for the little discs operad that operates on the Hochschild
cochains which reduces to the usual action on the
homology/cohomology level?}
\end{quest}

This question has an affirmative answer in many ways by picking a
suitable chain model for the little discs operad
\cite{Maxim,T,MS,Vor2,KS,MS2,BF}, see also \cite{MSS} for a review
of these constructions. We will provide a new, natural,
transparent and minimal positive answer to this question, by
giving an operation of $\CCcact$ on the Hochschild cochains.

There is a certain minimal set of operations necessary for the
proof of such a statement which is given by iterations of the
operations $\cup$ and $\circ_i$. These are, as we argue below, in
bijective correspondence with trees in $\wlbptree$. Our model for
the chains of the little discs operad $\CCcact$ has cells which
are exactly indexed by these trees hence the operation of the cellular
operad contains only the minimal number of operations and all operations
are non-zero. Furthermore, the top
dimensional cells which control the bracket constitute the universal
concatenation operad and hence yield operations for any operad,
see \S \ref{metaops}.

Finally, we will show that the differential which contracts arcs
of the cactus can be seen as a topological version of the
Hochschild differential. This makes our new topological solution
natural and minimal.

The main statement of this section is

\begin{thm}\label{delconj} Deligne's conjecture is true for the chain model of
the little discs operad provided by $\CCcact$, that is $CH^*(A,A)$
is a dg-algebra over $\CCcact$ lifting the Gerstenhaber algebra
structure.
\end{thm}

\begin{proof}
This follows from Theorem \ref{celldecomp} together
with either Proposition \ref{diffprop} or equivalently Proposition
\ref{foliageproof}.
\end{proof}

\begin{rmk} We will give two proofs
using the cell model $\CCcact$ for the chains of the little disc
operad by defining actions of its operadically isomorphic model
$\wlbptree$. First the trees act naturally on the Hochschild
complex by considering a tree as a flow chart of brace operations
and multiplications; the signs being fixed by one of the schemes
discussed below. Going one level deeper, instead of flow charts
for brace and multiplication operations we will use a so-called
foliage operator to reduce the brace operations to those of
insertion.
\end{rmk}

\subsection{The brace operations} The following operations appear
naturally when considering the iterations of Gerstenhaber's
operation $\circ$. They were first described by Getzler
\cite{Gbrace} and \cite{Kad} and are called brace operations. For
homogeneous $f,g_i$ of degrees $|f|$ and $|g_i|$, $N=|f| + \sum_i
|g_i| -n$

\begin{multline}
\label{bracedef}
  f\{g_1,\dots ,g_n\}(x_1,\dots,x_N) :=\\
\sum_{\footnotesize\begin{tabular}{c}
$1 \leq i_1 \leq \dots \leq i_n \leq |f|:$\\
 $i_j + |g_j|\leq i_{j+1}$\end{tabular}}\pm
f(x_1, \dots, x_{i_1-1},g_1(x_{i_1}, \dots, x_{i_1+|g_1|}), \dots, \\
\dots ,x_{i_n-1}, g_n(x_{i_n}, \dots, x_{i_n+|g_n|}), \dots , x_N)
\end{multline}
where the sign is the sign of the shuffle of the $g_j$ and $x_i$
which is determined by assigning shifted degrees to the $x_i$ and
$g_j$; namely the $x_i$ are considered to be of degree 1 and the
$g_j$ are considered to be of degree $|g_j|+1$.

Notice that $f\{g\}=f\circ g$.

\subsubsection{The suboperad generated by the brace operations}
It is well known (\cite{Gbrace, Kad, GV} see also \cite{MSS} for
the history of the brace operations and Deligne's conjecture) that
the set of concatenations of multiplications and brace operations
form a suboperad of the endomorphism operad of the Hochschild
complex. We will call it $\brace$.

The generators of this suboperad are in 1-1 correspondence with
elements of $\wlbptree$. Such a tree represents a flow chart. The
functions to be acted upon are to be inserted into the white
vertices. A black vertex signifies the multiplication of the
incoming entities. A white vertex represents the brace operation
of the element attached to that vertex outside the brace and the
elements obtained  by performing the operations associated to the
incoming branches inside the braces.

\begin{rmk}
Notice that in the flow chart of an expression of the type
$f\{(g_1),(g_2\cdot g_3 \cdot g_4\{h_1,h_2\})\}$ the symbols
``$\{$'' and ``$,$'' correspond to the white edges when reading
off the operation from a labelled tree $\t\in \wlbptree(n)$ in the
order $\prec^{\t}$.
\end{rmk}

\begin{prop}
\label{braceiso} The association of a flow chart is a non-$\Sigma$
operadic isomorphism between $\brace$ and $\wlbptree$ of operads
with a differential.
\end{prop}

\begin{proof}
The fact that the association of a flow chart is a bijection on
the generators was already mentioned. It is straightforward to
check the combinatorics of inserting a formula made up out of
braces and multiplications into a brace or a multiplication leads
to exactly the behavior described by our operad structure on the
trees $\wlbptree$. The checking of the compatibility of the signs
is tedious but straightforward. We would like to remark that for
the signs, in the brace formalisms the sums can be viewed as being
parameterized over a discretized simplex, in the sense that the
size of the gaps (number of variables between function insertions)
parameterize the summands in the formula (\ref{bracedef}) and the
sum of all the sizes of the gaps is fixed. This formalism also
yields the agreement of signs. The compatibility of the
differential can then be seen by a either a straightforward
calculation, a comparison with \cite{G,GV} or the above formalism
of discretized simplices.
\end{proof}

\subsection{Signs for $\brace$}

\begin{df}
We define an action of the symmetric group on $\brace$, by
considering the symbols ``$\{$'' and ``$,$'' to be each of degree
one.
\end{df}

\begin{prop}
With the above action of the symmetric group on $\brace$ the
isomorphism of \ref{braceiso} is an isomorphism of operads.
\end{prop}

\begin{proof}
This follows from the fact that the white edges correspond to the
symbols ``$\{$'' and ``$,$''. It can also be seen directly by
comparing to the formulas of \cite{G} and \cite{GV}.
\end{proof}

\subsection{The Operation of $\CCcact$ on $\Hom_{CH}$}
For $\O= \Hom_{CH}$ the endomorphism operad of $CH^*(A,A)$,
we define the map $\rho$ of equation
(\ref{opmaps}) to be given by the operadic extension of the maps
which send the tree $\t_n$ of figure \ref{nleaf} to the
non-intersecting brace operations. We let $\t_0$ act as the
identity and $\tau^b_n$ as multiplication.

As discussed above, in order to make signs match for the symmetric
group actions, we consider the action as described in \S
\ref{operadicaction}. To get complete agreement with the signs of
\cite{G}, we will have to consider the opposite orientation for
tensor factors in (\ref{tensordarstellung}) to that of $W$, i.e.\
we use the order $\overline{W}:= W_N \otimes \dots \otimes W_{1}$
for the tensor factors. To implement this change of sign we define
$\sign^W(\t)$ to be the sign obtained by shuffling the graded
tensor factors from the order $W$ to the order $\overline{W}$.
This basically means that in the orientation given by $W$ one
would regard the operations $\circ^{op}$ and $\cup^{op}$ on the
Hochschild complex, where $f\circ^{op}g= (-1)^{pq+p+1}g\circ f$
and $f\cup^{op}g:= (-1)^{pq} g\cup f$.

The action of the tree $\t_n$ is given by:
\begin{equation}
f \otimes L_{e_1} \otimes g_1  \otimes  \dots  \otimes  L_{e_n} \otimes  g_n
\mapsto (-1)^{\sign^W(\t)} f\{g_1,\dots ,g_n\}
\end{equation}
Here we used the notation $L_{e_1}$ to mean $L$ in the position of
$e_1$ in the tensor product of the form (\ref{tensordarstellung}).

The action of $\tau^b_n$ is given by
\begin{equation}
 g_1  \otimes  \dots  \otimes  g_n
\mapsto (-1)^{\sign^W(\t)} g_1 \cup \dots \cup g_n
\end{equation}

The operadic extension means that we read the tree  as a flow
chart: at each black vertex $|v|$ the operation $\t^b_{|v|}$
is performed and at each white vertex the operation $\t^w_{|v|}$
is performed. The $\Sn$ action is given by permutations and indeed
induces the right signs on the Hochschild complex as seen by straightforward
calculation.

\begin{prop}
\label{del1} The above procedure gives an operation of $\CCcact$
on $CH^*(A,A)$ and an operadic isomorphism between $\brace$ and
$\CCcact$.
\end{prop}

\begin{proof}
This follows directly from Theorem \ref{treescelliso}, the
previous paragraph and the considerations above.
\end{proof}

\subsubsection{The differential}
Denote the differential on $\CCcact$ by $\partial$ and the
differential of $CH^*(A,A)$ by $\delta$. On the space $W$ there is
a natural differential $\partial_W$ which is induced by
$\delta+\partial$. The differential on $W$ is induced by the
respective tree differential which equivalently collapses the arcs
of the spineless cactus or removes factors of $W$. Then the action
of $\CCcact$ on $\Hom_{CH}$ commutes with the differential in the
following sense.

\begin{prop}
\label{diffprop}
\begin{equation}
\rho \circ (\partial_W) = \delta\circ \rho
\end{equation}
The Hochschild cochains $CH^*(A,A)$ are a dg-algebra over
$\CCcact$ and there is an operadic isomorphism of the differential
operads $\brace$ and $\CCcact$.
\end{prop}

\begin{proof}
The verification of the compatibility of the grading and
differentials is a straightforward computation.
\end{proof}

\subsection{A second approach to the operation of $\CCcact$}
\label{maxcell} Another way to make $\CCcact$ or $\wlbptree$ act
is by using the foliage operator, see below \S \ref{foliage}. This
approach was first taken in \cite{KS}. It stresses the fact that a
function $f\in CH^n(A,A)$ is naturally depicted by $\t_n$. Notice
for instance the compatibility of the differentials. The tree for
$\del(f)$ is $\del(\t_n)$ where $\del$ is the differential on
trees with tails given in \S\ref{diffcheck}.

\subsubsection{Natural operations on $CH^*$ and their tree depiction}
\label{treepicture} Given elements of the Hochschild cochain
complex there are two types of natural operations which are
defined for them. Suppose $f_i$ is a homogeneous element, then it
is given by a function $f:A^{\otimes n} \rightarrow A$. Viewing
the cochains as functions, we have the operation of insertion. The
second type of operation comes from the fact that $A$ is an
associative algebra; therefore, for each collection $f_1,\dots,f_n
\in CH^*(A,A)$ we have the n! ways of multiplying them together.

We can encode the concatenation of these operations into a black
and white bipartite tree as follows: Suppose that we would like to
build a cochain by using insertion and multiplication on the
homogeneous cochains $f_1,\dots, f_n$. First we represent each
function $f_i$ as a white vertex with $|f_i|$ inputs and one
output with the cyclic order according to the inputs $1,\dots,
|f_i|$ of the function. For each insertion of a function into a
function we put a black vertex of arity one having as input edge
the output of the function to be inserted and as an output edge
the input of the function into which the insertion is being made.
For a multiplication of $k\geq 2$ functions we put a black vertex
whose inputs represent the factors in the order of their
multiplication. Finally, we add tails to the tree by putting a
black vertex at the end of each input edge which has not yet been
given a black vertex, and we decorate the tails by variables
$a_1,\dots,a_s$ according to their order in the total order of the
vertices of the rooted planted planar tree. It is clear that this
determines a black and white bipartite tree.

\subsubsection{Operations on Hochschild from trees with tails}
 A rooted planted planar bipartite black and white
tree whose tails are all black and decorated by variables
$a_1,\dots,a_s$ and whose white vertices are labelled by
homogeneous elements $f_v\in CH^{|v|}(A,A)$ determines an element
in $CH^{s}(A,A)$ by using the tree as a flow chart. This means the
operation of insertion for each black vertex of arity one and
multiplication for each black vertex of higher arity. Notice that,
since the algebra is associative, given an ordered set of elements
$a_i:i\in \{1,\dots,n\}$ there is a unique multiplication
$\prod_{i=1}^na_i$.

\begin{rmk} Using the above procedures, the possible ways
to compose $k$ homogeneous elements of
$CH^*(A,A)$ using insertion and cup product are bijectively
enumerated by black and white bipartite planar rooted planted
trees with tails and $k$ white vertices labelled by $k$ functions
whose degree is equal to the arity of the vertex.
\end{rmk}

\begin{nota}
We will fix $A$ and use the short hand notation $CH:= CH^*(A,A)$.
For an element $f\in CH$, we write $f^{(d)}$ for its homogeneous
component of degree $d$.
\end{nota}

\begin{df}
For $\t\in \wlbptreet (n)$ and $f_1, \dots, f_n\in CH$ we let
$\t(f_1, \dots, f_n)$ be the operation obtained in the above fashion
by decorating the vertex $v_i$ with label $i$ with the
homogeneous component of $f_i^{(|v_i|)}$.
Notice that the result is zero if any of the homogeneous
components  $f_i^{(|v_i|)}$
vanishes.
\end{df}

\begin{rmk}
Up to the signs which are discussed below this gives an operation
of $\CCcact$ on the Hochschild complex.
\end{rmk}

\subsection{The operation of $\wlbptree$}

In order to define the operation
we need foliation operators in the botanical sense, i.e.\ operators
that add leaves.
To avoid confusion with the mathematical term ``foliation'', we choose to
abuse the English language and call these operations ``foliage'' operators.

\begin{df}
\label{foliage} Let $l_n$ be the tree in $\wlbptreet$ with one
white vertex labelled by $v$ and $n$ tails as depicted in figure
\ref{nleaf}. The foliage operator $F:\wlptree \rightarrow
\wlptreet$ is defined by the following equation
$$
F(\t):= \sum_{n\in \mathbb {N}}l_n \circ_v \t
$$
\end{df}

\begin{rmk}
Notice that the right hand side is infinite, but  $\wlptreet$ is
graded by  the number of leaves, and $F(\t)$ is finite for a fixed
number of black leaves so that the definition does not pose any
problems. Furthermore, one could let $F$ take values in $
\wlptreet[[t]]$ where $t$ keeps track of the number of tails which
would make the grading explicit.

Also notice that $F:\wlbptree\rightarrow \wlbptreet$ and
$F:\wptree\rightarrow \wpttree$.
\end{rmk}

\begin{df}
\label{foliageopdef} For a tree $\t\in \wlbptree(n)$  with $n$
white vertices we define a map $\op(\t)\in \mathrm
{Hom}(CH^{\otimes n}, CH) = \Hom_{CH}(n)$ by
$$
op(\t)(f_1, \dots, f_n):= \pm (F(\t))(f_1, \dots, f_n)
$$
here the signs are as discussed in \S \ref{genfix} and the right
hand side is well defined since it only has finitely many non-zero
terms.
\end{df}

\begin{rmk}
The considerations of this section naturally lead to brace
operations in a far more general setting. This is explained in
detail in \S \ref{metaops}.
\end{rmk}

\subsubsection{Differential on trees with tails}
\label{diffcheck}

\begin{df}
 For a tree $\t$ with tails in $\bwbptree$
and vertex $v \in V_b\setminus \{v_{root}\}$ we define $\t^+_v$ to
be the b/w tree obtained by adding a black vertex $b+$ and an edge
$e^+:=\{b+,v\}$, if $|v|\neq 0$ and if $|v|= 0$, the tree obtained
by adding two vertices $b+$ and $b_{st}$ and two edges
$e+=\{b+,v\}$ and $e_{st}=\{b_{st},v\}$ to $\t$.

We call a linear order $\prec'$ on $\t^+_v$ compatible with the
order $\prec$ on $\t$ if a) $e+\prec' e_{st}$ if applicable and b)
the order induced on $\t$ by $\prec'$ by contracting $e^+$ and
$e_{st}$ (if applicable) coincides with $\prec$. We define
$E_{w-int}$ to be the white internal edges, i.e. white edges which
are not leaves and set $E_{b-angle}: E(\t) \setminus E_{w-int}$.
For a linear order $\prec'$ on $\t^+_v$
$$\sign(\prec'):=(-1)^{|\{e|e
\in E_{b-angle}, e\prec'v \}|}$$ and set
$$\del_v(\t):=\sum_{\text{compatible}\prec'}\sign(\prec')(\t^+_v,\prec').$$
We recall that tail edges are considered to be black. Finally, we
define

\begin{equation}
\label{taildiff}\del(\t):=\sum_{v \in V_b\setminus
\{v_{root}\}}\del_v(\t)
\end{equation}
\end{df}

\begin{rmk}
In the tree depiction for the operations of inserting and cup
product \S \ref{treepicture}, the tree differential amounts to
inserting $\cup$ products into the ``slots'' represented by black
vertices. Using this interpretation and the tree notation for the
calculations of \cite{G,GV} it is straightforward to check that
the tree differential (\ref{taildiff}) defined above agrees with
the differential on the Hochschild co-chains.
\end{rmk}

\begin{prop}
\label{foliageproof}
 The Hochschild cochains are a $dg$-algebra over the operad
$\wlbptree$.
\end{prop}

\begin{proof}
The properties of an operad follow in a straightforward way from
the definition \ref{foliageopdef} and the definition of the
differentials.
\end{proof}

\subsection{The top dimensional cells of spineless cacti and pre-Lie operad}

 We denote the  top-dimensional cells of $\CCncactn$ by
$CC_n^{top}(n)$. These cells again form an operad and they are
indexed by trees with black vertices of arity one.  Furthermore,
the symmetric combinations of these cells which are the image of
$\flrtree$ under the embedding $\cppin$  form a sub-operad, see
Proposition \ref{cppinop}. For one choice of orientation we make
the signs explicit in the next Lemma.

\begin{lem}
\label{symtreeop} In the orientation $\blab$ for the
top-dimensional cells for $\t \in \flrtree(n)$ and $\t'\in
\flrtree(m)$
$$ \cppin(\t)\circ_i\cppin(\t')=\pm\cppin(\t\circ_i\t')
$$
where the sign is $(-1)^{(i-1)(n-1)}$ if the root vertex of $\t$
has a label which is less than $i$ and  $(-1)^{i(n-1)}$ if the
root vertex of $\t$ has a label which is bigger than $i$.
\end{lem}

\begin{proof}
First notice that by Proposition \ref{cppinop}, we have that the
right hand side contains the terms indexed by the trees appearing
in the embedding. Furthermore, notice that under the choice of
signs induced by the orientation $\blab$ the signs do not depend
on the particular structure of the tree and are only dictated by
the labels. In the composition the labels are such that the $n$
labels of the second tree are permuted to the $i$-th label of the
first tree. Therefore there is a universal sign which is given by
$(-1)^{(i-1)(n-1)}$ if the root vertex of $\t$ has a label which
is less than $i$ and by $(-1)^{i(n-1)}$ if the root vertex of $\t$
has a label which is bigger than $i$.
\end{proof}

\begin{df} Let $\mathcal{G}Pl$ be the quadratic operad in the
category $Vect_{\mathbb{Z}}$ obtained as the quotient of a free
operad $\mathcal{F}$ generated by the regular representations of
$\mathbb{S}_2$ by the quadratic relations defining a graded
pre-Lie algebra. Let $\mathcal{R}$ be the ideal  generated by the
graded $\mathbb{S}_3$ submodule
 generated by the relation
$$r=(x_1x_2)x_3-x_1(x_2x_3)
- (-1)^{|x_2||x_1|} ((x_1x_3)x_2 -x_1(x_3x_2)).$$ Then
$\mathcal{G}Pl=\mathcal{F}/\mathcal{I}$. Here $\mathcal {F}$ and
$\mathcal{R}$ are considered to be graded by assigning the degree
$n-1$ to $\mathcal{F}(n)$.
\end{df}

This operad is the operad for graded pre-Lie algebras in the sense
that any algebra over this operad is a graded pre-Lie algebra and
vice-versa any graded pre-Lie algebra is an algebra over this
operad:

\begin{thm}
\label{preliecells} The operad $CC_n^{top}(n)^{\mathbb {S}}\otimes
k$ is isomorphic to the operad $\mathcal{G}Pl$ for graded pre-Lie
algebras. Furthermore the shifted operad $(CC_n^{top}\otimes
(L^*)^{\otimes E_w})^{\mathbb {S}}(n)\otimes k$ is isomorphic to
the operad $\mathcal{P}l$ for pre-Lie algebras.

The analogous statements also hold over $\mathbb{Z}$.
\end{thm}

\begin{proof}
In view of Lemma \ref{symtreeop} and Proposition \ref{cppinop},
the second statement follows from the operadic isomorphism of
$\flrtree$ and the pre-Lie operad $\mathcal{P}l$ of \cite{CL}.
This also proves the first statement up to signs.  The matching of
the signs is guaranteed by the shift, see Definition
\ref{shifted}. The fact that the relation $r$ holds and generates
the respective ideal was verified in the presentation of
Gerstenhaber structure for spineless cacti \cite{cact} by a
translation from the explicitly given relation on the chains of
the arc operad \cite{KLP}.
\end{proof}

\begin{cor}
\label{preliecor} The  pre-Lie algebra of $\Sn$-coinvariants
$\bigoplus_n ((CC_n^{top}\otimes (L^*)^{\otimes E_w})^{\mathbb
{S}}(n))_{\Sn}$ is isomorphic to the free pre-Lie algebra in one
generator.

Likewise the graded pre-Lie algebra of  $\Sn$-coinvariants
$\bigoplus_n(CC_n^{top})^{\mathbb {S}}(n))_{\Sn}$ is isomorphic to the graded
free pre-Lie algebra in one generator
\end{cor}

\begin{proof}
The first statement follows from \cite{CL} and thus so does the
second up to signs. These are guaranteed to conform by the
shifting procedure and Theorem \ref{preliecells}.
\end{proof}

\section{Structures on Operads and Meta-Operads}
\label{metaops} In this paragraph, we discuss how Deligne's
conjecture and the Gerstenhaber structure for the Hochschild
complex can in fact be generalized to structures on operads. This
helps to explain some choices of signs and explains the naturality
of the construction of insertion operads which gives a special
role to spineless cacti as their topological incarnation as well
as to $\Arc$ as a natural generalization.

This analysis also enables us surprisingly to relate spineless
cacti to the renormalization Hopf algebra of Connes and Kreimer
\cite{CK}, see the next paragraph.

\subsection{The universal concatenations}
Given any operad there are certain universal operations, i.e.\
maps of the operad to itself. We will first ignore possible signs
and comment on them later.

Any operad possesses the operations given by its structure maps
$$
\circ_i: O(m)\otimes O(n) \rightarrow O(m+n-1)
$$
Any concatenations of these maps will also yield operations on the
operad.  All possible concatenations of the structure map are
 described by their flow charts.
 These charts are in turn given by trees
$\t \in \wpttree$ as we now explain.

Any concatenation of $k$ objects $op_i\in O(n_i)$ using the
structural maps $\circ_i$ will be given by
\begin{equation}
\label{gentreeops} (\dots((op_{\sigma(1)}\circ_{i_1}
\op_{\sigma(2)})\circ_{i_2}op_{\sigma(3)} \circ_{i_3} )\dots
)\circ_{i_{k-1}} op_{\sigma(k)}
\end{equation}
 with $i_1<i_2,\dots
<i_k$ and some permutation $\sigma \in \mathbb{S}_k$.

Let $\t \in \wpttree(k)$  such that  $ |v_i|=n_i, i\in 1, \dots k$ and set $m=\sum_{i=1}^k n_i-k-1$
then there is an operation
\begin{equation}
\label{treeops}
\circ(\t): (O(n_1) \otimes \dots \otimes O(n_k))\rightarrow O(m)
\end{equation}
which is obtained as follows. First label the vertex $v_i$ by
$op_{n_i}\in O(n_i)$. Let $\sigma \in \mathbb{S}_k$ be the
permutation which is defined by the order of the tree, i.e.\
$v_{\sigma(1)}\prec^{\t} \dots \prec^{\t}v_{\sigma(i)}\prec^{\t}
v_{\sigma(i+1)}\prec^{\t} \dots \prec^{\t} v_{\sigma (k)}$.

Now starting  at the bottom of the tree and going along the
outside edge path, we read off the operation
\begin{equation}
 \circ(\t)(op_1\otimes \dots \otimes
op_k):=(\dots((op_{\sigma(1)}\circ_{i_1}
\op_{\sigma(2)})\circ_{i_2}op_{\sigma(3)} \circ_{i_3} )\dots
)\circ_{i_{k-1}} op_{\sigma(k)}
\end{equation}
where the $i_k$ are given by $i_k=|\{v\in E_b(\t):v\prec^{\t}
v_{\sigma(k)}\}|$; recall that the root is considered a black
vertex.

Vice-versa, given a concatenation as above we can successively
build up the corresponding tree.

\begin{rmk} Notice, we might have white leaves, which
allows one to consider operads with a $0$ component --- such as
$CH^*$. In general lifting the restriction on the $n_i$, we define
the operations $\circ(\t)$ to be zero if  $|v_i| \neq n_i$.
\end{rmk}

\subsubsection{The insertion partial meta-operad}
Algebraically we can concatenate the operations $\circ(\t)$ by
substituting
$$op_i=\circ(\t')(op'_1\otimes \dots \otimes op'_l) \text{ for some
$\t'\in \wpttree(l)$} $$ into (\ref{gentreeops}). In general, for
$\t \in \wpttree(k)$, let $v_i:=\lb_{\t}(i)$ be the $i$-th white
vertex of $\t$, and let $ \t'\in \wpttree(l)$ with $|v_i|=|V_{
leaf}(\t')|$, we define the tree $\t \circ_i \t'$ to be the tree
of the concatenated operation
\begin{multline}
\label{partialops}
\circ(\t\circ_i \t')(op_1\otimes \dots \otimes op_{k+l-1}):= \\
\circ(\t)(\op_1\otimes \dots \otimes op_{i-1} \otimes \circ(\t')(op_{i}
\otimes \dots \otimes op_{i+l-1})\otimes
op_{i+l} \otimes \dots \otimes op_{l+k-1})
\end{multline}

\begin{lem}
The multiplication maps $\circ_i$ (\ref{partialops}) together with
the permutation action on the labels imbue  $\wpttree$ with the
structure of a  partial operad.\footnote{A partial operad is the
notion obtained from an operad by requiring that the composition
maps $\circ_i$ are only defined on a subobject of the $O(n)$.
These compositions should be equivariant with respect to the $\Sn$
action and satisfy the associativity axioms if it is possible to
compose them.}
\end{lem}
\begin{proof}
Straightforward.
\end{proof}

\begin{rmk}
The partial concatenations $\circ_i$ insert a tree with $k$-tails
into the vertex $v_i$ if $|v_i|=k$, by connecting the incoming
edges of $v_i$ to the tail vertices in the linear order at $v_i$
and then contracting the images of the tail edges.
\end{rmk}

\begin{df}
\label{admitssum} We will fix that for $\O$ in $\Set$ the direct
sum which we again denote by $\O$  is given by the free Abelian
group generated by $\O$ which we consider to be graded by the
arity of the operations $\op \in \O$ minus one. If the operad $\O$
is in $\Chain$ we can take the direct sum of the components as
$\mathbb{Z}$-modules. In the case of an operad $\O$ in the
category $\Vect$ we consider its direct sum to be the direct sum
over $k$ of its components. In all these cases, we say
$\O=\{O(n)\}$ admits a direct sum and write $\O=\bigoplus_{n\in
\mathbb{N}}O(n)$. We also consider $\O$ to be graded by
$\mathbb{N}$ with the degree of $O(n)$ being $n-1$.
\end{df}

\begin{rmk}
Consider an operad which admits a direct sum. Let $\O$
be its direct sum, then we obtain a  map of partial operads.
$$
\wpttree \rightarrow Hom(\O,\O)
$$

In this sense one can say that $\wpttree$
is the universal concatenation partial meta-operad.
\end{rmk}

\subsection{The pre-Lie Structure of an Operad}
\label{prelieop}

 In an operad which admits a direct sum,
 one can define the analog of the $\circ$ product of \cite{G} and the iterated
brace operation (cf.\ \cite{Gbrace,GV,GV2}).

\begin{df}
Given any operad $\O$ in $\Set$, $\Chain$ or $\Vect$, we define
the following map
\begin{eqnarray}
O(m) \otimes O(n) &\rightarrow& O(m+n-1)\\
op_m\otimes op_n &\mapsto& \sum_{i=1}^m (-1)^{(i-1)(n+1)} op_m \circ_i op_n
\end{eqnarray}
This extends to a map on $\O=\bigoplus_{n \in \mathbb{N}} O(n)$
\begin{equation}
 \circ:\O\otimes \O\rightarrow\O
\end{equation}
which we call the $\circ$ product.

 We call the map which is defined in the same fashion as $\circ$, but
 with the omission of the signs $(-1)^{(i-1)(n+1)}$  the ungraded $\circ$
product.
\end{df}

\begin{prop}
The product $\circ$ defines the structure of a graded pre-Lie
algebra on $O:= \bigoplus_{i\in \mathbb{N}}O(n)$. Omitting the
sign $(-1)^{(i-1)(n+1)}$ in the sum yields the structure of a
non-graded pre-Lie algebra.
\end{prop}

\begin{proof}
The proof is a straightforward calculation which is analogous to
Gerstenhaber's original calculation  \cite{G}. We do not wish to
rewrite the proof here, but in graphical notation
 the proof follows from figure
\ref{prelieprop} below. Without signs this notation is related to
the one that can be found in \cite{CL}, for the case with signs
the interpretation of the resulting trees is discussed in
\S\ref{genfix} and Proposition \ref{cppinop}.
\end{proof}

See also \cite{MSn} for another variant of this fact.

\begin{prop}
If an operad admits a direct sum then its direct sum is an algebra
over the symmetric top dimensional chains of the little disc
operad of the chain model provided by $\CCcact$ as well as over
the shifted chains $(CC_n^{top})^{\mathbb {S}}\otimes
(L^*)^{\otimes E_w}$.
\end{prop}

\begin{proof}
We have shown that the direct sum of an operad which admits such a
sum has the structure of a pre-Lie algebra and a graded pre-Lie
algebra so that the statement follows from \ref{preliecells}.
\end{proof}

\subsection{The insertion operad}

The interesting property of the operation $\circ$ is that it
effectively removes the dependence on the number of inputs of the
factors. Using this logic systematically, we obtain a universal
insertion operad.

\subsubsection{Actions of $\wptree$}
Given an operad which admits a direct sum, we can also define
other operations similar to $\circ_i$. The summands of these
operations are brace operations, which are in natural
correspondence with $\wptree$. In fact these operations all appear
in the iterations of $\circ$. They are given by inserting the
operations into each other according to the scheme of the tree. In
other words, we will show that every operad is a brace algebra.

Essentially, at this level of abstraction, we would not like to { a priori}
specify the number of leaves, i.e.\ inputs and
degrees of the operations, so we  have to consider trees
with all possible decorations by leaves. For this, we can use
the foliage operators of Definition \ref{foliage}.

Recall that there is an operation of $\wpttree$ on homogeneous
elements of $\O$ of the right degree. We extend this operation to
all of $\O$ by extending linearly and setting to zero expressions
which do not satisfy degree condition. Where the degree condition
corresponding to a given $\t$ applied to homogenous elements
$\op_k$ is that $op_k \in O(|v_k|)$.

Given $\t$ in $\wptree(n)$, we then define the operation

\begin{equation}
\label{foilops} \circ(\t)(op_1 \otimes \dots\otimes op_n):=
F(\t)(op_1\otimes \dots \otimes op_n)
\end{equation}

Notice that, although $F(\t)$ is an infinite linear combination,
for given $op_1, \dots, op_n$ the expression on the right hand
side is finite.

\begin{figure}
\epsfbox{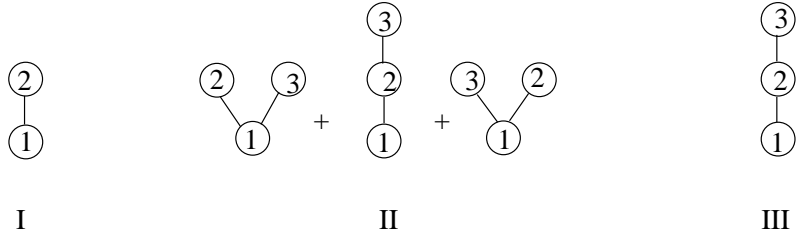} \caption{\label{prelieprop}I. $op_1 \circ
op_2$ II. $(op_1\circ op_2)\circ op_3$ and III. $op_1\circ
(op_2\circ op_3)$}
\end{figure}

 Examples of this are given in
figure \ref{prelieprop}. Here the first tree yields the operation
$op_1 \circ op_2$, i.e the insertion (at every place) of $op_2$
into $op_1$. Iterating this insertion we obtain expression II
which shows that inserting $op_3$ into $op_1 \circ op_2$ gives
rise to three topological types: inserting $op_3$ in front of
$op_2$, into $op_2$ and behind $op_2$. In the opposite iteration
one just inserts $op_2 \circ op_3$ into $op_1$ which gives a
linear insertion of $op_2$ into $op_1$ and $op_3$ into $op_2$.
From the figure (up to signs) one can read off the symmetry in the
entries 2 and 3 of the associator. The signs are fixed by the
considerations of \S \ref{genfix}.

\begin{rmk}
\label{insproduct}

Again, inserting $op_i=\circ(\t')(op'_1\otimes \dots\otimes
op'_k)$ we obtain operad maps with respect to the symmetric group
actions on the labels
$$
\circ_i:\wpttree(k)\otimes \wptree(l) \rightarrow \wptree(k+l-1)
$$
by demanding that for all $op_1,\dots, op_{k+l-1}$
\begin{multline*}
\circ(\t\circ_i \t')(op_1 \otimes \dots \otimes op_{k+l-1})=\\
\circ(\t)(op_1\otimes \dots \otimes op_{i-1} \otimes
\circ(\t')(op_{i} \otimes \dots \otimes op_{i+l-1})\otimes
op_{i+l} \otimes \dots \otimes op_{l+k-1})
\end{multline*}
Or in other words thinking about $F$ as a formal power-series,
e.g.\ in $\wpttree[[t]]$ where the powers of $t$ keep track of the number of tails,
$$
F(\t \circ_i \t')=F(\t) \circ_i F(\t')
$$

Thus we have an insertion operad structure on $\wptree$.
\end{rmk}

\begin{rmk}
We can induce a pre-Lie operation on $\wpttree$ via the $\circ$
operation.
\begin{equation}
F(\t_1*\t_2):= F(\t_1) \circ F(\t_2)
\end{equation}
\end{rmk}

\begin{thm}
Any chain operad which admits a direct sum is an algebra over the
operad $\wptree$ with the insertion products.
\end{thm}

\begin{proof}
Immediate from the preceding.
\end{proof}

\subsection{Operad algebras and a generalized Deligne conjecture}

\begin{df}
We define an operad algebra to be an operad $\O$ which admits a
direct sum together with an element $\cup\in O(2)$ which is
associative, i.e.\ if $a\cup b:=(-1)^{|a|}(\cup\circ_1
a)\circ_{|a+1|} b$ then $(a\cup b)\cup c = a \cup (b\cup c)$.
Recall that $|a|=n-1$ if $a\in O(n)$.
\end{df}

This definition is essentially equivalent to the definition of an
operad with a multiplication of \cite{GV2}.

\begin{df}
For $op\in O(m), op'_i\in O(n_i)$, we define the generalized brace
operations
\begin{multline}
\label{genbrace}
  op\{op'_1,\dots ,op'_n\} :=
\sum_{\footnotesize\begin{tabular}{c}
$1 \leq i_1 \leq \dots \leq i_n \leq m:$\\
 $i_j + |op'_j+1|\leq i_{j+1}$\end{tabular}}\pm
(\cdots((op\, \circ_{i_1} op'_1)\circ_{i_2}op'_2)\circ_{i_3}
\cdots )\circ_{i_n} op'_n
\end{multline}
where the sign is defined to be the same one as in equation.
(\ref{bracedef})
\end{df}

\begin{lem}
\label{opactionlem}
 There is an operadic action of $\wlbptree$ of any operad
algebra.
\end{lem}

\begin{proof} We can view the bipartite tree as a flow chart. For
the white vertices, we use the brace operations (\ref{genbrace})
and for a black vertex with $n$ incoming edges, we use the $n-1$
times iterated operation $\cup$. Notice that the order in which we
perform these operations does not matter, since we took $\cup$ to
be associative.
\end{proof}

\begin{defprop}
\label{gerstdiff} Generalizing Gerstenhaber's \cite{G} definition
to an operad algebra, we define a differential on the direct sum
by $\del f = f\circ \cup - (-1)^{|f|} \cup \circ f$.
\end{defprop}

\begin{proof} The fact that this is a differential follows from the
calculations of \cite{G}.
\end{proof}

\begin{rmk}
The analogous tree differential is given in \S\ref{diffcheck}.
Replacing functions by elements of the operad the compatibility
follows for the more general setup of Definition \ref{gerstdiff}.
\end{rmk}

\begin{thm}\label{gendel}
The generalized Deligne conjecture holds. I.e.\  the direct sum of
any operad algebra which admits a direct sum  is an algebra over
the chains of the little discs operad in the sense that it is an
algebra over the dg-operad $\CCcact$.
\end{thm}

\begin{proof}
By the preceding Lemma \ref{opactionlem}, we have an operadic
action of $\wlbptree$ and thus an action of the chains $\CCcact$
which is a chain model for the little discs operad. The
compatibility of the differentials follows directly from their
definitions by a straightforward calculation as remarked above.
\end{proof}

\section{The Hopf algebra of Connes and Kreimer and spineless Cacti}

\subsection{The Hopf algebra of an operad}
We have seen in \S\ref{prelieop} that any operad that admits a
direct sum gives rise to a pre-Lie algebra. Now by the defining
properties for a pre-Lie algebra  the commutator of its product
yields a Lie algebra or in the graded case an odd Lie algebra.

By the above considerations, we can thus naturally associate a
pre-Lie algebra, a Lie algebra, and a Hopf algebra to each operad
that admits a direct sum.

\begin{df} Given an operad $\{O(n)\}$ which admits a direct
sum $\O=\bigoplus_n O(n)$, we define its pre-Lie algebra
$\mathcal{PL}(\O)$  to be the pre-Lie algebra $(\O, \circ)$ with
$\circ$ as defined in \S\ref{prelieop}, its Lie algebra
$\mathcal{L}ie(\O)$ to be the Lie algebra $(\O, [ \, ,\, ])$using
the Lie bracket $[a,b]:=a\circ b-b\circ a$ and its odd Lie algebra
$\mathcal{L}ie^{\Z/2\Z}(\O)$ to be the odd Lie-algebra $(\O, \{ \,
, \, \})$ where $\{ \, , \, \}$ is defined as usual via $ \{a \,
,\, b\}:= a\circ b - (-1)^{(|a|+1)(|b|+1)}b\circ a$ for $a\in
O(|a|)$ and $b\in O(|b|)$. Finally, the Hopf algebra of an operad
$\mathcal{H}opf(\O)$ is defined to be $U^*(\mathcal{L}ie(\O))$,
i.e.\ the dual of the universal enveloping algebra of its Lie
algebra.
\end{df}

As it turns out these objects or their $\Sn$-coinvariants are of
interest. In particular, we obtain chain models in terms of chains
of spineless cacti or even moduli space for some well known
algebras.

\subsection{Connes-Kreimer's
Hopf Algebra as the Hopf algebra of an operad} In \cite{CK} a Hopf
algebra based on trees was defined to explain the procedure of
renormalization in terms of the antipode of a Hopf algebra. This
Hopf algebra was described directly, but also as the dual to the
universal enveloping algebra of certain Lie algebra which was
identified in \cite{CL} as the Lie algebra associated to the free
pre-Lie algebra in one generator.

\begin{df}
By the $\mathbb{S}_n$ coinvariants of an operad which admits a
direct sum, we mean $\bigoplus_{n\in \mathbb{N}}(O(n))_{\Sn}$. We
write $\O_{\mathbb{S}}$ for this sum. Here $\bigoplus_{n\in
\mathbb{N}}(O(n))_{\Sn}$ is the shorthand notation explained in \S
\ref{admitssum}.
\end{df}

In our notation we can rephrase the results of \cite{CK,CL} as

\begin{thm}
\label{clthm} The renormalization Hopf algebra of Connes and
Kreimer $H_{CK}$ is isomorphic to the Hopf algebra of
$\Sn$-coinvariants of $\H(\flrtree)$. This Hopf algebra is also
isomorphic to the $\Sn$-coinvariants of $\H(\mathcal{P}l)$.
\end{thm}

\subsection{A chain interpretation of $H_{CK}$} \label{ckcell}

As shown in Theorem \ref{preliecells} there is a cell and thus a
topological interpretation of the pre-Lie operad and the graded
pre-Lie operad inside $\cact^1$ and thus inside the $\Arc$ operad.
In this interpretation $H_{CK}$ is also the Hopf algebra of the
coinvariants of the shifted chain operad
$CC^{top}_*(\cact)^{\mathbb S}\otimes (L^*)^{\otimes E_w}$. Recall
that $L^*$ is a free $\mathbb{Z}$ module generated by an element
$l$ of degree $-1$.

\begin{prop}\label{ckcor}
$H_{CK}$ is isomorphic to the Hopf algebra of $\Sn$ coinvariants
of the sub-operad of top-dimensional symmetric combinations of
shifted cells $((\CCcact)^{top})^{\mathbb{S}}\otimes
(L^*)^{\otimes E_w}$ of the shifted cellular chain operad of
normalized spineless cacti $CCcact\otimes (L^*)^{\otimes E_w}$
$$
H_{CK}\simeq (\mathcal{H}opf((\CCcact)^{top})^{\mathbb{S}}\otimes
(L^*)^{\otimes E_w})_{\mathbb{S}}.
$$
\end{prop}
\begin{proof}
Immediate from Theorem \ref{clthm} and Corollary \ref{preliecor}.
\end{proof}

 It is interesting to note that also the G- and BV-structures
which are given by spineless cacti and cacti \cite{cact} are
inside the symmetric (graded symmetric) combinations as well.

\begin{df}
We define the planar Connes Kreimer Hopf algebra $H^{pl}_{CK}$ to
be $(\mathcal{H}opf(CC_*^{top}\otimes (L^*)^{\otimes E_w}))_{\mathbb{S}}$. This
is the straightforward generalization from rooted to planar
planted trees.
\end{df}

\begin{rmk}
We can also consider the graded Hopf algebra corresponding to the
$\Sn$-coinvariants of graded top dimensional cells
$U^*(\mathcal{L}ie^{\Z/2\Z}(\O))_{\mathbb{S}}$.
\end{rmk}

\begin{rmk}
The examples above should be relevant in the context of multiple
polylogarithms. We expect to obtain other interesting examples of
such Hopf algebras by applying the above constructions to other
operads based on trees.
\end{rmk}

\subsection{Comments on Operads and $H_{CK}$}

We have shown that any operad is an algebra over the operad
$\flrtree$ in a natural way and thus the Hopf algebra $H_{CK}$
naturally appears in any context involving operads, such as
Deligne's conjecture. We have furthermore shown that there is a
topological incarnation of the insertion product, which is based
on surfaces. In this setting, we have constructed a chain
representation of the algebra $H_{CK}$. This links the algebra
$H_{CK}$ and its underlying bracket for instance to string
topology and in positive characteristic to Dyer-Lashof-Cohen
operations \cite{Tour}.

\section{Further developments and generalizations}

\subsection{Generalizations}
\subsubsection{The $A_{\infty}$-case}
All the statements made above in the context of algebras hold as
well for $A_{\infty}$-algebras. For this one has to change to an
equivalent cell model for normalized spineless cacti. A sketch of
the procedure is given below. A detailed account will be given in
\cite{KSch}.

For $\t \in \wlptree$
let
\begin{equation*}
\label{cellinfty} \cell(\t) := \prod_{v\in V_w} W_{|v|+1} \times
\prod _{v\in V_b} K_{|v|}
\end{equation*}
where $W_{|v|+1} $ is the $|v|$--dimensional cyclohedron and
$K_{|v|}$ is the $|v|-2$--dimensional Stasheff polytope a.k.a.\
associahedron.

Now one has to show a compatibility:
\begin{lem}
The tree differential $\del^{\infty}_{tree}$ of \cite{KS} and the
natural differential $\del_{poly}$ for the cyclohedra and
associahedra are compatible in the following sense. The cells
$\cell(\t')$ appearing in the sum $\cell(\del^{\infty}_{tree} \t)$
are in 1-1 correspondence to the products appearing in
$\del_{poly}\cell(\t)$ and furthermore the signs agree.
\end{lem}

\begin{df}
Let $Cell(\wlptree)$ be the CW complex glued from the cells
$\cell(\t)$ using $\del^{\infty}_{tree}$.
\end{df}

\begin{prop}
$CC_*(Cell(\wlptree))$ is equivalent to $\CCcact$ and hence a cell
model for the little discs operad.
\end{prop}
\begin{proof}
Just contract the associahedra simultaneously blowing
 down the cyclohedra to simplices. The compatibility of
this operation follows from the properties of the map $\forgetainfty$.
\end{proof}

Using flow charts as in \cite{KS} where now a black vertex $b$
corresponds to the $|b|$-th of the $A_{\infty}$-multiplication
$\mu_{|b|}:A^{\otimes |b|}\rightarrow A$ one obtains:

\begin{thm}\cite{KSch}
Deligne's conjecture holds in the $A_{\infty}$ case over
$\mathbb{Z}$ for the chain model $CC_*(Cell(\wlptree))$.
Furthermore so does its generalization to an operad with an
$A_{\infty}$-multiplication.
\end{thm}

Here an operad with an $A_{\infty}$ multiplication is an operad
$\O$ with a collection of multiplications elements $m_n\in O(n)$
which satisfy the equations for $A_{\infty}$ multiplications.

\begin{rmk}
Using compatible realizations of cyclohedra and associahedra given as particular faces, one can also construct a topological quasi-operad \cite{KSch}. An interesting open problem is the definition of
a suboperad of $\Arc$ which is the natural  topological operad model for these cells.
\end{rmk}

\subsubsection{The cyclic case}

The results and methods of this paper have been extended to the
cyclic case in \cite{cyclic} in which it is proven that the
Hochschild cochains of a Frobenius algebra carry an action of a
chain model of the framed little discs operad. For this one needs
other techniques than presented in this paper, notably operadic
correlation functions.

\subsection{Further developments}
Building on the techniques of \cite{cyclic} one can extend the
action to include all marked ribbon graphs and all $\Z/2\Z$
decorated ribbon graphs \cite{hoch1,hoch2}. The latter generalizes
all known graph actions on the Hochschild co-chains such as the
chains framed little discs and Sullivan Chord diagrams. These
results yield an action of pseudo-chains of the decorated moduli
space. There remains one caveat, since the ribbon graphs only give
pseudo-cells, the dg-structure at the moment still needs some
clarification.

In another new development, by taking up the philosophy of the
present paper, we are able to realize all little $k$-cubes operads
inside the arc formalism by using stabilization with respect to
the genus operator of \cite{KLP}. Using the formalism developed in
this paper we then obtain a solution to Kontsevich's
generalization of Deligne's conjecture to $d$-algebras. Moreover
in the stable limit, we obtain an $E_{\infty}$ suboperad
\cite{loop}. This implies that the stabilization of the $\Arc$
operad yields an infinite loop space spectrum. We would like to
point out, that this stabilization is different from the usual
stabilization. For this one would use the gluing on of a fixed
element of a particular type of the $\Arc$ operad which
is not considered a
stabilization in the above sense. Nevertheless, the proposed
results of \cite{loop} could be regarded as a combinatorial/chain
incarnation of the seminal results of Tillmann \cite{Til,Til2} and
Madsen-Tillmann \cite{MT}.

\section*{Acknowledgments}
I would like to thank the Max-Planck Institute for Mathematics and
the IHES for their hospitality. The visits to these Institutes
were instrumental in the carrying out of the presented research;
especially the discussions with Maxim Kontsevich. Furthermore, I
would like to thank Ieke Moerdijk for interesting discussions. My
thanks also goes to Wolfgang L\"uck and Peter Teichner for their
hospitality and  interest. In fact, the first steps for our
solution to Deligne's conjecture were taken during a visit of the
SFB in M\"unster and the final step occurred to me during a visit
to San Diego in response to questions of Peter Teichner. Finally,
I wish to thank Jim Stasheff for very enlightening discussions and
encouragement as well as Murray Gerstenhaber for pointing out the
similarities between his constructions and those of \cite{KLP}.
This observation was the initial spark for this work on Deligne's
conjecture. We also wish to thank the referee for his careful
reading and helpful suggestions.

The author also acknowledges partial financial support from the NSF
under grant \#DMS-0070681.

\renewcommand{\theequation}{A-\arabic{equation}}
\renewcommand{\thesection}{A}
\setcounter{equation}{0}  
\setcounter{subsection}{0}
\section*{Appendix: Graph-theoretic Details}  
\subsection{Dual graph constructions}
\subsubsection{The dual bipartite b/w graph}
Fix a  marked spineless treelike ribbon graph $(\G,f_0)$. We
define  a new bipartite b/w graph $\t(\G)$ as follows. Set
$$V_b(\t(\G)):=V(\G),\; V_w(\t(\G)):=\{\mbox{cycles of
$\G$}\}\setminus \{c_0\},\; F(\t(\G)):=F(\G)$$ and specify $f_0$
as the distinguished flag. Now set
$$
\delta_{\t(\G)}(f):=\begin{cases}\delta_{\G}(f)& \mbox{if } f\in
c_0\\
c_i & \mbox{if } f\in c_i \neq c_0\end{cases}, \quad
\imath_{\t(\G)}(f):=\begin{cases} N_{\G}(\imath_{\G}(f)) &
\mbox{if } f\in c_0 \\
\imath_{\G}(N^{-1}_{\G}(f)) &\mbox{if }f\notin c_0.\end{cases}$$
 We fix the ribbon graph structure
by declaring the cyclic order at the white vertices to be the
conjugate order of the cycle they represent. That is $f$ is the
predecessor of $f'$ in the cyclic order at $v=c_i$ if $f$ is the
successor of $f'$ in the cycle $c_i$. For the black vertices, we
take the induced order from the identifications
$V_b(\t(\G))=V(\G)$ and $F(\t)=F(\G)$. This means in particular
that
$$N_{\t(\G)}(f):=\begin{cases}
\imath_{\G}(f)&\mbox{if }f\in c_0\\N_{\G}(\imath_{\G}(f))&\mbox{if
} f\notin c_0\end{cases}$$
Notice that this graph only has one cycle and hence has genus $0$,
so it is a planar tree. This first statement can be seen directly
or as follows, we have that the number of vertices of $\t$ is
equal to the number of vertices of $\G$ plus the number of cycles
of $\G$ minus 1. The number of flags of $\t$ is the number of
flags of $\G$ which is 2 times the number of edges of $\G$. So we
get that $2-2g(\t)=(|V(\G)|+\#\mbox{cycles of
$\G$}-1)-E(\G)+\#\mbox{cycles of $\t$}=1+\#\mbox{cycles of
$\t$}\leq2$.

\subsubsection{Ribbon graph defined by a b/w planar planted bipartite tree}
The ribbon graph defined by a tree is given as follows: Given a
b/w bipartite planar planted tree $\t$, we define a new graph
$\G(\t)$ by setting
$$F(\G(\t)):=F(\t),\quad
V(\G(\t)):=V_b(\t)$$ and specifying $f_0$ as the distinguished
flag. Furthermore, we define
$$\delta_{\G(\t)}(f):=\begin{cases}\delta_{\t}(f)& \mbox{if }\delta_{\t}(f)\in
V_b(\t)\\
\delta_{\t}(\imath_{\t}(f))& \mbox{if }\delta_{\t}(f)\in
V_w(\t)\end{cases}, \imath_{\G(\t)}(f):=\begin{cases}
N_{\t}&\mbox{if }\delta_{\t}(f)\in V_b(\t)\\N_{\t}^{-1}& \mbox{if
}\delta_{\t}(f)\in V_w(\t).\end{cases}$$ We fix the ribbon graph
structure by defining the cyclic order at a vertex as follows: if
$\delta(f)\in V_b$ then the successor of $f$ is the flag
$\imath_{\t}(f)$, the successor of this flag is defined to be the
successor of $f$ in the cyclic order at $v=\delta(f)$ of $\t$ and
so on. This means that if $\delta(f)\in V_w$ then the successor of
$f$ is $N_{\t}(f)$. This implies that
$$N_{\G(\t)}(f):=\begin{cases}
N^2_{\t}(f)&\mbox{if }\delta(f)\in
V_b\\
\imath_{\t}(N^{-1}_{\t}(f))&\mbox{if } \delta(f) \in
V_w\end{cases}$$

Now $\G$ has a cycle which runs through all the edges and the
other cycles of $\G$ are in 1-1 correspondence with the white
vertices. The cycle that contains $f_0$ contains every second flag
in the unique cycle of the tree, that is all the flags $f$ with
$\delta(f)\in V_b$.  On the other hand if $\delta(f')=v\in V_w$
then all the flags in its cycle are incident to the same vertex
and in effect $\imath_{\t}( N^{-1}_{\t}(f)$ is just the
predecessor of $f$ in the cyclic order of its vertex. Thus these
cycles are in 1--1 correspondence with the sets $F(v)$ for $v\in
V_w$. This means that $2-2g(\G)=|V_b(\t)|-|E(\t)|+|V_w(\t)|+1=2$,
so that the genus of $\G$ is zero. The cycle of $f_0$ contains all
flags $f$ with $\delta_{\t}(f)\in V_b$ and for such a flag
$\delta_{\t}(\imath_{\G}(f))\in V_w$ and vice versa, so that
either $f$ or $\imath_{\G}(f)$ lie in the $N_{\G}$ cycle of $f_0$
and indeed the graph $\G$ is a marked spineless treelike ribbon
graph.

\begin{lem}The dual graph is a duality transformation that is
$\G(\t(\G))=\G$ and $\t(\G(\t))=\t$.\end{lem}

\begin{proof}
  On the level of flags  this is
clear. For the vertices, this is also clear for the first order of
iteration, in the second order the equality of the set of vertices
follows from the observation explained above that the cycles of
$\G(\t)$ which are not the distinguished cycle $c_0$ are in 1-1
correspondence to the white vertices. It remains to check the
compatibility of the maps $\delta,\imath$ and $N$ which amounts to
plugging in the definitions.  Here are some examples:
$\delta_{\G(\t(\G))}(f)$ with $f\in
c_0$:$\delta_{\G(\t(\G))}(f)=\delta_{\t(\G)}(f)=\delta_{\G}(f)$.
If $f\in c_i\neq c_0$ then
$\delta_{\G(\t(\G))}(f)=\delta_{\t(\G)}(\imath_{\t(\G)}(f))=\delta_{\t(\G)}(N_{\G}(\imath_{\G}(f)))=\delta_{\G}(N_{\G}(\imath_{G}(f)))=\delta_{\G}(f)$.
As another example consider $N_{\t(\G(\t))}(f)$ for $f$ with
$\delta_{\t}(f)\in V_b$, we get
$N_{\t(\G(\t))}(f)=\imath_{\G(\t)}(f)=N_{\G}(f)$ and if
$\delta_{\t}(f)\in V_w$ then
$N_{\t(\G(\t))}(f)=\imath_{\G(\t)}(f)=N_{\G(\t)}(\imath_{\G(\t)}(f))=N_{\G(\t)}(N^{-1}_{\t}(f))=N^2_{\t}(N^{-1}(f))=N_{\t}(f)$.
Writing out the other calculations is tedious but straightforward.
\end{proof}
\subsection{Gluing spineless and normalized spineless cacti: the graph version}
\subsubsection{$S^1$--graphs}
To formulate the gluing in a purely graph theoretic way, we first need a new definition. An {\em $S^1$-graph}  is a metric ribbon graph of genus $0$ together with  a  distinguished flag  $f_0$, such that
\begin{itemize}
\item [i)] if $c_0$ is the cycle of $f_0$ then
for each flag either $ f\in c_0$ or $\imath(f)\in c_0$
\item[ii)]  $\forall v: |F(v)|=2$
\end{itemize}

Such a graph has exactly two cycles that have the same length, which we call the radius of the $S^1$-graph.

An $S^1$--graph is equivalent to the data of a circle embedded
into the plane with several marked points. To explain this, let
$S^1_r:=\{(x,y)\in {\mathbb R}^2|x^2+y^2=r\}$. We will use the
natural coordinate $\theta$ on $S^1$:
$x=r\cos(2\pi\theta),y=r\sin(2\pi\theta)$. Let $0=
\theta_0<\theta_1 < \dots< \theta_n< 1$ be $n+1$--points on $S^1$.
Then $(S^1,(\theta_i))$ defines an $S^1$ graph by letting the
points be the vertices, the arcs  the edges and the flags are the
half edges. To be very explicit each flag is an ordered pair
$(\theta_i,\theta_j)$ with $|i-j|=1 \equiv (n-1)$ and
$\imath(v,w)=(w,v)$. We denote the edges by
$\{\theta_i,\theta_j\}$ with the same restriction. Since
$|F(v)|=2$ there is only one choice of cyclic order at each
vertex. We let $f_0=(\theta_0,\theta_1)$ and let $c_0$ be the
cycle $f_0$ that contains $f_0$. Finally, we set
$\mu(\{\theta_0,\theta_n\}):=1-\theta_n$ and
$\mu(\{\theta_i,\theta_{i+1}\})=\theta_{i+1}-\theta_{i}$ for
$i=0,\dots n-1$. Conversely, any $S^1$ graph gives rise to the
data above. For this we embed the realization of the CW complex
which is an $S^1$ into the plane as a circle $S^1_r$ of radius
equal to the length of $c_0$, such that $c_0$ runs
counterclockwise and that the vertex of $f_0$ is the point
$\theta_0=0$.  Each vertex then corresponds to a point $\theta_i$.

\subsubsection{Gluing $S^1$--graphs}
\label{sgluing} We define a gluing for two $S^1$--graphs of the
same radius. Given $S$ and $S'$, we scale $S'$, so that the
lengths of the radii of $S$ and $S'$ agree and then let $S\circ
S'$ be the $S^1$ graph defined by taking the union of marked
points. That is if $S=(S^1_r,\{ \theta_i:i\in I\})$ and
$S'=(S^1_{r},\{\theta'_j:j\in J\})$ then $S\circ
S'=(S^1_r,\{\theta_i:i\in I\}\cup \{\theta'_j:j\in J\})$.

We identify the set of vertices of the glued graphs with the union
of the vertices of the two graphs: $V(S\circ S')=V(S)\cup V(S')$,
we can also naturally identify $F(S\circ S')$ with $F(S)\cup
F(S')$.

We extend the gluing of $S^1$--graphs to glueings of an
$S^1$-subgraph with an $S^1$--graph.\footnote{A subgraph is  a
subset of vertices and a subset of flags closed under
$\imath_{\G}$ and $\delta_{\G}$.} Let $(\G,\mu)$ be a metric
ribbon graph and $S$ be a subgraph which is an $S^1$--graph with
the induced ribbon graph structure. Also fix an $S^1$--graph $S'$,
we define $(\G,\mu) \circ_{S}S'$ by replacing the subgraph $S$ by
$S\circ S'$. This is well defined, since we can identify the
vertices of $S$ with vertices of $S\circ S'$. The result is again
a metric ribbon graph.

This operation simply inserts new vertices into the subgraph, new vertices will have valence $2$, so that it is naturally a ribbon graph.

\subsubsection{Lobes as $S^1$--graphs}
Every lobe gives rise to an $S^1$--graph by considering the
subgraph of vertices and edges of the cycle. We define the marked
flag $f_i$ of the $S^1$--graph to be the first flag $f$ of the
cycle $c_0$ such that $\imath(f)\in c_i$.

The fact that a cycle corresponding to a lobe is actually an
$S^1$--subgraph, that is that each vertex $v$  has  $|F(v)|=2$
when considering the subgraph, can most easily be seen from the
dual tree. Here the statement corresponds to the fact that any two
edges incident to a white vertex only have this vertex in common,
which is true since the dual graph is a tree.

\subsubsection{Chord diagrams}
There is an $S^1$--graph covering each spineless cactus as
follows. For a spineless cactus $c=(\G,f_0,\mu)$  consider the
$S^1$--graph obtained by going around the cycle $c_0$ and
considering this an $S^1$--graph. The set of flags is defined to
be the same as that of $\G$ and $\imath$ and $\mu$ are also
defined identically in both graphs, but let the set of vertices of
the $S^1$-graph be the set $F(\G)$ and fix the map $\delta=id$.
The marked flag is again $f_0$. We call this graph $S(c)$.

There is an equivalence relation on the vertices of $S(c)$ given
by $v=\delta_{S(c)}(f)\sim v'=\delta_{S(c)}(f')$ if
$\delta(f)=\delta(f')$ in $c$. We call $(S(c),\sim)$ the chord
diagram of $c$. We wish to point out that as graphs $\G=S(c)/\sim$
and that $(S(c)/\sim,f_0,\mu)$ is $c$ if we give the $S(c)/\sim$
the ribbon structure of $\G$ that is the structure induced by
first lifting $N_{\G}$ as a function on flags and then letting it
descend to the quotient graph.

\subsubsection{The glueings for spineless cacti}

Fix two spineless cacti  $c=(\G,f_0,\mu)$ and
$c'=(\G',f_0',\mu')$, let $S$ be the $S^1$--subgraph of $c$ and
represent $c'$ as $(S(c'),\sim)$. Fix a cycle $c_i$ of $c$ and let
$r_i$ be the radius of $c_i$ and $R$ be the radius of $c'_0$, the
outside circle of $c'$.

We define the underlying metric  graph of $c\circ_i c'$ to be the
graph given by
\begin{equation}[(\G,\mu) \circ_{c_i}\frac{r_i}{R}S(c')]/\sim
\end{equation}
where $\frac{r_i}{R}$ is the scaling action and $\sim$ is the
equivalence relation given above extended to the graph $(\G,\mu)
\circ_{c_i}\frac{r_i}{R}S(c')$ by identifying the vertices of
$\frac{r_i}{R}S(c')$ as a subset of the vertices of the glued
graph. We mark this graph by the image of the flag $f_0$ of $c$.
The ribbon structure is defined as follows. For all vertices of
$\G$ which do not lie on $c_i$, we keep the order. For the other
vertices, we fix the following linear order. Given $v$ let
$v_1,v_2,\dots, v_n$ be its pre-images in $(\G,\mu)
\circ_{c_i}\frac{r_i}{R}S(c')$ linearly ordered according to the
distinguished cycle $c_0'$ of $S(c')$ and its distinguished flag.
Each of the sets $F(v_i)$ is also linearly ordered by using its
cyclic order and declaring the first flag to be the first flag of
$c'_0$ which is incident to the vertex.
 Now $F(v)=\coprod_{i=1,\dots,n} F(v_i)$ and we enumerate the
flags by the induced linear order starting to enumerate $F(v_1)$
in its linear order, continuing with $F(v_2)$ and finishing with
$F(v_n)$ again in their linear orders. This gives a linear order
on $F(v)$ and we fix the cyclic order of $F(v)$ to be the unique
cyclic order which is compatible with that linear order.

 It is now straightforward
to check that the resulting data defines a spineless cactus. With
the exception of $c_i$, all the cycles which were lobes  descend
unaltered as subsets of flags. The Euler--characteristic
computation then shows that the genus is zero and that there is
exactly one more cycle. This is the cycle of the flag $f_0$ and it
passes through all the edges. It follows that the number of lobes
is sub-additive.

\subsubsection{The glueings for normalized spineless cacti}
Keeping the notations of the paragraph above,  now let $c$ and
$c'$ be normalized spineless cacti and $c_i$ a fixed cycle of $\G$
with radius $r_i$. Then set  $\tilde\mu(e)=\mu(e)$ if $e$ is not
in $c_i$ that is none of the two flags is in $c_i$. If one of the
flags of $e$ is in $c_i$ then $\tilde\mu(e)=\frac{R}{r_i}(e)$.
This causes the radius of $c_i$ considered as a cycle of $\tilde
c$ to be $R$, the radius of the outside circle of $c'$. We then
define the gluing for normalized spineless cacti by setting the
underlying metric graph to be
\begin{equation}
[(\G,\tilde \mu)\circ_{c_i}S(c')]/\sim.
\end{equation}
Again we mark the resulting graph by the image of the flag $f_0$
and induce a ribbon graph structure as above. As above it follows
 that the result is a normalized spineless cactus and the gluing is sub-additive in the
 number of lobes.

\end{document}